\documentclass{gtart_a}
\pdfoutput=1
\usepackage{pinlabel}

\title[Asymptotic geometry of the mapping class group]{Asymptotic geometry of the mapping class group\\and Teichm\"uller space}

\author{Jason A Behrstock}
\givenname{Jason A}
\surname{Behrstock}
\address{Department of Mathematics\\Barnard College\\\newline
Columbia University\\New York, NY 10027\\USA}
\email{jason@math.columbia.edu}
\urladdr{}

\volumenumber{10}
\issuenumber{}
\publicationyear{2006}
\papernumber{38}
\lognumber{0641}
\startpage{1523}
\endpage{1578}

\doi{}
\MR{}
\Zbl{}

\arxivreference{math.GT/0502367}


\received{15 August 2005}
\revised{15 July 2006}
\accepted{26 September 2006}
\proposed{Benson Farb}
\seconded{Walter Neumann, Rob Kirby}
\accepted{}
\published{21 October 2006}
\publishedonline{21 October 2006}
\corresponding{}
\editor{CPR}
\version{}

\AtBeginDocument{\let\bar\wbar
\def\C{{\mathcal C}}\def\notin{\not\in}}
\makeop{Graph}

\numberwithin{equation}{section}

\makeatletter
\def\cnewtheorem#1[#2]#3{\newtheorem{#1}{#3}[section]
\expandafter\let\csname c@#1\endcsname\c@thm}

\newtheorem{thm}{Theorem}[section]
\cnewtheorem{prop}[thm]{Proposition}
\cnewtheorem{cor}[thm]{Corollary} 
\cnewtheorem{conj}[thm]{Conjecture}
\cnewtheorem{qn}[thm]{Question} 
\cnewtheorem{lem}[thm]{Lemma}
\cnewtheorem{cvn}[thm]{Convention}
\cnewtheorem{fthm}[thm]{Folk Theorem}
\newtheorem*{thm*}{Theorem}

\theoremstyle{definition} 
\cnewtheorem{defn}[thm]{Definition}
\cnewtheorem{nota}[thm]{Notation}
\cnewtheorem{rmk}[thm]{Remark}
\cnewtheorem{exa}[thm]{Example}
\cnewtheorem{observation}[thm]{Observation}

\makeatother  

\makeautorefname{defn}{Definition}

\def\H{{\mathbb H}}
\def\N{{\mathbb N}}

\def\G{{\mathcal G}}

\def\MCG{\mathcal {MCG}}

\def\M{{\mathcal M}}

\def\I{{\mathcal I}}
\def\P{{\mathcal P}}

\def\<{\langle}
\def\>{\rangle}

\def\cone{\mathrm{Cone}_{\omega}}
\def\coneomega{\mathrm{Cone}_{\omega}}
\def\con{\mathrm{Cone}_{\omega}}
\def\coneM{\M_{\omega}}

\def\ulim{\lim_\omega}
\def\ulimi{\lim_\omega}
\def\ulimif{\lim_{\omega}\frac{1}{d_{i}}}

\def\seq{\mathrm{Seq}}
\def\diam{\mathrm{diam}}

\def\phihat{\widehat{\Phi}}

\newcommand{\base}{\operatorname{base}}

\def\I{{\mathbf I}}
\def\T{{\mathbf T}}
\newcommand{\fsub}{\mathrel{\scriptstyle\searrow}}
\newcommand{\bsub}{\mathrel{\scriptstyle\swarrow}}
\newcommand{\fsubd}{\mathrel{{\scriptstyle\searrow}\kern-1ex^d\kern0.5ex}}
\newcommand{\bsubd}{\mathrel{{\scriptstyle\swarrow}\kern-1.6ex^d\kern0.8ex}}
\newcommand{\fsubeq}{\mathrel{\raise-.7ex\hbox{$\overset{\searrow}{=}$}}}
\newcommand{\bsubeq}{\mathrel{\raise-.7ex\hbox{$\overset{\swarrow}{=}$}}}

\begin{document}

\begin{asciiabstract}
In this work, we study the asymptotic geometry of the mapping class
group and Teichmuller space. We introduce tools for analyzing the
geometry of `projection' maps from these spaces to curve complexes of
subsurfaces; from this we obtain information concerning the topology
of their asymptotic cones. We deduce several applications of this
analysis. One of which is that the asymptotic cone of the mapping
class group of any surface is tree-graded in the sense of Drutu and
Sapir; this tree-grading has several consequences including answering
a question of Drutu and Sapir concerning relatively hyperbolic
groups. Another application is a generalization of the result of Brock
and Farb that for low complexity surfaces Teichmueller space, with the
Weil-Petersson metric, is delta-hyperbolic. Although for higher
complexity surfaces these spaces are not delta-hyperbolic, we
establish the presence of previously unknown negative curvature
phenomena in the mapping class group and Teichmueller space for
arbitrary surfaces.
\end{asciiabstract}

\begin{abstract} 
In this work, we study the asymptotic geometry of the mapping class group
and Teichmuller space. We introduce tools for analyzing  the geometry of 
``projection'' maps from these spaces to curve complexes of subsurfaces;
from this we obtain information concerning the topology of their
asymptotic cones. We deduce several applications of this analysis. One 
of which is that the asymptotic cone of the mapping class 
group of any surface is tree-graded in the
sense of Dru\c{t}u and Sapir; this tree-grading has several
consequences including answering a question of 
Dru\c{t}u and Sapir concerning relatively hyperbolic groups. Another
application is a  generalization of the result of Brock and Farb 
that for low complexity
surfaces Teichm\"{u}ller space, with the Weil--Petersson metric, is
$\delta$--hyperbolic. Although for higher complexity surfaces these
spaces are 
not $\delta$--hyperbolic, we establish the  presence of previously 
unknown negative curvature phenomena in the mapping class group and 
Teichm\"{u}ller space for arbitrary surfaces. 
\end{abstract}
\maketitle

\section{Introduction}

\subsection{Overview of main results}

In studying the geometry and dynamics of a surface $S$, some of the
most important mathematical objects are its mapping class group 
$\MCG(S)$, its Teichm\"{u}ller space, and the complex of curves. 
We prove new results about the asymptotic geometry of these 
objects through a study of their asymptotic cones.

Recently, Dru\c{t}u and Sapir 
showed that any group which is (strongly) relatively hyperbolic has 
cut-points in its asymptotic cone \cite{DrutuSapir:TreeGraded}. They 
asked whether this property of having asymptotic cut-points 
gives a  characterization of such groups. We prove the following.

\medskip
{\bf\fullref{cutpointstw}}\qua 
{\sl If $S$ is any surface and $\cone(\MCG(S))$ is any asymptotic cone,   
then every point of  $\cone(\MCG(S))$ is a global cut-point. 

Similarly one obtains global cut-points in each asymptotic cone of  
Teichm\"{u}ller space with the Weil--Petersson metric.}
\medskip

This result answers the question of \cite{DrutuSapir:TreeGraded} in
the negative for, as shown by Anderson, Aramayona and Shackleton
\cite{AASh:RelHypMCG} and independently by the author with Dru\c{t}u
and Mosher \cite{BehrstockDrutuMosher:thick}, the group $\MCG(S)$ is
not hyperbolic relative to any family of finitely generated subgroups.
This establishes the mapping class groups as the first known examples
of groups which are not relatively hyperbolic, but whose asymptotic
cones are tree-graded.  See \cite{BehrstockDrutuMosher:thick} for
other applications of \fullref{cutpointstw} and its
\fullref{MCGtreegraded}.

We  obtain the following more specific result on the asymptotic 
geometry of the mapping class group.

\medskip
{\bf \fullref{Fzertree}}\qua
{\sl In any surface, $S$, and any asymptotic cone, $\cone(\MCG(S))$, 
there exists an isometrically embedded $\R$--tree, 
$F_{0}\subset\cone(\MCG(S))$ with the following properties:
\begin{itemize}
    \item $F_{0}$ contains infinitely many distinct 
    infinite geodesics arising as asymptotic cones 
    of axes of pseudo-Anosovs;

    \item for each pair of points on $F_{0}$, there exists a unique 
    topological path in\break $\cone(\MCG(S))$ connecting them and 
    this path is contained in $F_{0}$.
\end{itemize}

A similar result holds replacing the mapping class group by    
Teichm\"{u}ller space with the Weil--Petersson metric.}
\medskip

One application of \fullref{Fzertree} is that it shows that  
the directions of pseudo-Anosov axes are negatively curved, ie, 
in the direction of each pseudo-Anosov 
the divergence function is superlinear. This readily implies 
that these directions are quasi-geodesically stable, ie, for the 
quasi-geodesic axis given by taking powers of 
a pseudo-Anosov, there is a  bound $M(K,C)$ for which any 
$(K,C)$--quasigeodesic with endpoints on this path must be completely 
contained within the $M$ neighborhood of the axis. 
We also establish that the Weil--Petersson metric on Teichm\"{u}ller 
space is negatively curved in the directions of geodesics all of whose 
subsurface projections are uniformly bounded. 
Thus although for higher complexity surfaces the Weil--Petersson metric 
is not $\delta$--hyperbolic, this result shows that there are 
directions in which one gets the useful consequences of
$\delta$--hyperbolicity. (Note that such a result does not follow from
the negative sectional curvature of the Weil--Petersson metric since
these curvatures are not bounded away from 0.) 
Results on stability of certain families of geodesics in the 
Teichm\"{u}ller metric was proved by Minsky \cite{Minsky:quasiprojections}, 
thus our results complete this analogy between these three spaces.
(After completing this result, we were informed by Minsky 
that in current joint work with Brock and Masur they will give an
alternate proof of the stability of a similar family of 
Weil--Petersson geodesics.)

When the complexity of the surface is small enough, we recover the 
following result which was first established by Brock and Farb 
\cite{BrockFarb:curvature}. 
(After this work was completed another alternate proof of this
result  appeared using interesting techniques involving the  
geometry of $CAT(0)$ spaces, see  Aramayona \cite{Aramayona:thesis}.)

\medskip
{\bf\fullref{pantshyperbolic}}\qua
{\sl With the Weil--Petersson metric, 
   the Teichm\"{u}ller spaces for the surfaces $S_{1,2}$ and $S_{0,5}$
   are each $\delta$--hyperbolic.} 
\medskip

This result implies that any subset of an 
asymptotic cone of Teichm\"{u}ller space of  $S_{1,2}$ or $S_{0,5}$ 
contains a point which disconnects it. By contrast, the asymptotic cone 
of the Teichm\"{u}ller space of any surface of higher complexity contains 
bi-Lipschitz embeddings of $\R^{n}$ with $n\geq 2$, these subsets can 
not be disconnected by any point. Nonetheless, \fullref{cutpointstw} 
implies that each point is a global cut-point; this is one of several 
ways to consider \fullref{cutpointstw} as a generalization of 
\fullref{pantshyperbolic}.

The original motivation for this work was to study quasiflats in the
mapping class group. In particular, it had been conjectured by
Brock--Farb that the mapping class group admits quasiflats only up to
the maximal dimension for which it contains abelian subgroups (by work
of Birman, Lubotzky and McCarthy \cite{BirmanLubotzkyMcCarthy} the
latter number is known to be $3g+p-3$).  Using \fullref{Fzertree} as
the base case for an inductive argument, the author in joint work with
Minsky \cite{BehrstockMinsky:rankconj} have proven this ``Rank
Conjecture'' (using different techniques Hamenst\"{a}dt
\cite{Hamenstadt:TT3} has also obtained a proof of this result).

\subsection{Methods}
The above results are produced by analyzing the \emph{complex of
curves}, a simplicial complex encoding intersection patterns of curves
on a surface. This complex was first introduced by Harvey and has
since been used to prove many deep results about the mapping class
group (see Harer \cite{Harer:vcd}, Harvey \cite{Harvey:Boundary},
Hatcher \cite{Hatcher:Pants}, and Ivanov \cite{Ivanov:mcg}). More
recently Masur and Minsky proved $\delta$--hyperbolicity of the complex
of curves \cite{MasurMinsky:complex1}, and then used this to establish
new methods for studying the mapping class group
\cite{MasurMinsky:complex2}. The complex of curves has also found many
applications to the study of 3--manifolds; in particular it has proved
crucial in the recent proof of the Ending Lamination Conjecture (see
Brock, Canary and Minsky \cite{Minsky:ELC1,BrockCanaryMinsky:ELC2}).
For other applications see: Behrstock and Margalit
\cite{BehrstockMargalit:superinj}, Bestvina and Fujiwara
\cite{BestvinaFujiwara:boundedcohom}, Bowditch
\cite{Bowditch:TightGeod}, Hempel \cite{Hempel:CurveComplexes}, Irmak
\cite{Irmak:SI1} and Margalit \cite{Margalit:AutPants}.

For each subsurface $X\subset S$, we consider  a ``projection'' map 
$$\pi_{X}\co\MCG(S)\to\C(X)$$\noindent from the 
mapping class group of $S$ into the 
curve complex of $X$, which is closely
related to the map taking a curve in $S$ to its intersection with $X$. 
The vertices of the curve complex of a surface, $\C(X)$, are homotopy 
classes of essential, non-peripheral  closed curves on $X$. The 
distance between two distinct homotopy classes of curves, 
$\mu,\nu\in\C(X)$, is a measure of how complicated the intersection of 
representatives of these curves must be. For example as long as $X$ 
has genus larger than one:  $d_{X}(\mu,\nu)=1$ if and only if $\mu$ 
and $\nu$ can be realized disjointly on $X$ whereas 
$d_{X}(\mu,\nu)\geq 3$ if and only if $\mu$ and $\nu$ fill $X$ 
(\emph{fill} means that for any representatives of $\mu$ and $\nu$ the 
set $X\setminus(\mu\cup\nu)$ consists of disks and once-punctured disks).

The following tool provides the starting point for our analysis.

\medskip
{\bf \fullref{projest}}\qua(Projection estimates)\qua
    {\sl Let $Y$ and $Z$ be two overlapping subsurfaces of $S$. Then for 
    any $\mu\in\MCG(S)$:
    $$d_{\C(Y)}(\partial Z, \mu)>M  \implies d_{\C(Z)}(\partial Y, 
    \mu)\leq M.$$
    Where $M$ depends only on the topological type of the surface $S$.}
\medskip

Combined with work of Masur--Minsky \cite{MasurMinsky:complex2} this 
result is useful in studying the metric geometry of mapping class 
groups and Teichm\"{u}ller space via their subsurface projections. 
(\fullref{projest} can be generalized to study projections of the 
mapping class group into a product of curve complexes. Doing so 
yields an interesting picture of certain aspects of the
asymptotic geometry of the mapping class group, this direction is 
developed further by the author in \cite{Behrstock:thesis}.)
Our study of the metric geometry of mapping class 
groups and Teichm\"{u}ller space is done via an analysis of their 
asymptotic cones. We define a projection map to send elements of the 
mapping class group of a surface $S$ to the mapping class group of a 
subsurface $Z\subset S$. This projection map has several nice 
properties, for instance it sends points to uniformly bounded 
diameter subsets (\fullref{diamboundmarkingproj}), and is used in  
defining an important new object introduced in this paper: the 
\emph{strongly bounded subset}, $F_{0}$, of the asymptotic cone. 
This subset is part of the \emph{transversal tree} of $\cone\MCG(S)$ 
and plays a large role in the  Behrstock--Minsky 
proof of the Rank Conjecture \cite{BehrstockMinsky:rankconj}.

\subsection{Outline of subsequent sections}

In \fullref{sectionbackgroundgeometry} we discuss the 
 tools of asymptotic geometry 
which are fundamental to the constructions and results in the rest of
this work. All the results 
in this section are well known and thus we are as much setting up 
notation as we are reminding the reader of which tools to keep at their 
fingertips for the remainder of this work. 

Surfaces are the main characters throughout this work. In
\fullref{chaptersurfacesbackground} we discuss the marking complex,
which is our quasi-isometric 
model of choice for dealing with the mapping class group. We also 
remind the reader of various facts we will need concerning 
the mapping class group, the complex of curves, and Teichm\"{u}ller space.

The original material in this work begins with the results in 
\fullref{chapterProjectionEstimates}. Here we calculate estimates 
on the image of the ``projection'' map from the marking 
complex into the product of the  curve complexes of the constituent 
subsurfaces. In this section we produce our starting point, 
\fullref{projest} (Projection estimates), 
which is a technical result we use in our analysis of 
the mapping class group. In 
\fullref{geomprojest} (Projection estimates; geometric 
version) we provide a geometric interpretation of 
\fullref{projest}.

The final three sections apply the Projection estimates Theorem and 
related techniques to various situations.
In the first of these we (re)prove hyperbolicity 
of the Weil--Petersson metric on Teichm\"{u}ller space and 
of the mapping class group in some low 
complexity cases, Theorems~\ref{pantshyperbolic} and 
\ref{mcghyperbolic}. The techniques introduced in 
this section provide a good warm-up for those used to prove the key
step in its generalization, \fullref{Fzertree}.

In \fullref{sectiontransversaltree} we generalize the hyperbolicity of low 
complexity mapping class groups and Teichm\"{u}ller spaces to 
deduce the existence of negatively curved directions in the cases of higher 
complexity surfaces.
The main result of this section is \fullref{Fzertree} which
identifies some new aspects of the geometry of $\cone(\MCG)$ which we 
exploit in the following section.

We end with \fullref{sectiontreegraded} where 
we relate our earlier work to questions about relative hyperbolicity
in the mapping class group. Here we answer some questions about 
mapping class groups and tree-graded spaces and raise some new question.

\subsection*{Acknowledgments}
The first five sections of this paper were part of my doctoral thesis
\cite{Behrstock:thesis}.  It is a pleasure to thank my thesis advisor
Y~Minsky for many useful discussions and for sharing with me his
enthusiasm for mathematics.  I would also like to thank L~Mosher for
sharing with me his knowledge of $\MCG$ and for his thorough reading
and comments on this work. I extend my thanks to many others whose
input and suggestions influenced this work---especially C~Leininger
and D~Margalit who have generously spent hours talking with me about
various aspects of the mapping class group.  I am grateful to the
Barnard and Columbia Math Departments for providing an ideal work
environment.  Finally, I would like to especially thank my wife,
Kashi, for her love and support.

\section{Background I: Large scale geometry}\label{sectionbackgroundgeometry}

\subsection{Geometric Group Theory}

There are several equivalent definitions for hyperbolicity; throughout
this paper we will have in mind Rips' {\em thin triangle}
characterization:
\begin{defn} A geodesic metric space $X$ is called
\emph{$\delta$--hyperbolic} if there exists a constant $\delta$ so that 
for each triple $a,b,c\in X$ and each choice of 
geodesics $[a,b], [b,c],$ and $[a,c]$ one has that $[a,b]$ is contained 
in a $\delta$--neighborhood of the union of $[a,c]$ and $[b,c]$.
\end{defn}

Note that we write $[a,b]$ to denote a geodesic between $a$ and $b$ 
even though this choice need not be unique.

\begin{defn}
Let $(X,d_X)$ and $(Y,d_Y)$ be metric spaces. 
A map $\phi\co X\to Y$ is called a \emph{$(K,C)$--quasi-isometric
 embedding}
 if there exist constants $K\geq 1$ and $C\geq  0$ 
 such that for all $a,b\in X$
$$
 \frac{1}{K}d_X(a,b)-C \leq d_Y(\phi(a),\phi(b))\leq K d_X(a,b) +C. 
$$
It is a \emph{$(K,C)$--quasi-isometry}\index{quasi-isometric} 
if additionally, $\phi$ has the 
property that every point of $Y$ lies in the $C$--neighborhood
of $\phi(X)$.

When there exists a $(K,C)$--quasi-isometry $\phi\co X\to Y$ between two 
metric spaces $X$ and $Y,$ we say that they are {\em 
$(K,C)$--quasi-isometric} and denote this $X\sim_{K,C}Y$; 
When the choice of constants is not important (as will often be the 
case), we simply say $X$ and $Y$ are
{\em quasi-isometric}, denoted $X\sim Y$. 

\end{defn}

For geodesic spaces hyperbolicity is a quasi-isometry invariant; note
that the hyperbolicity constant may change (see Ghys and de la Harpe
\cite{GhysHarpe:afterGromov}).

In order to apply this definition to groups, we recall the natural left 
invariant metric on a group $G$ with finite generating set 
$S$. The \emph{word metric on $G$ relative to 
the generating set $S$}\index{word metric} is given 
by $d_{G,S}(a,b)=|a^{-1} b|$ where 
$|a^{-1} b|$ is the smallest number of letters needed to represent the group 
element $a^{-1} b$ in terms of letters from $S$ and their inverses.
Although this metric 
depends on the choice of generating set, for finitely generated 
groups, its quasi-isometry type does not. Thus since hyperbolicity is 
a quasi-isometry invariant, the notion of a group being {\em hyperbolic} 
is well defined.

A basic non-hyperbolic example  where we will use the notion of
quasi-isometries is $\Z^{n}$ (with $n>1$),  
which we (easily) observe is quasi-isometric to $\R^{n}$.
In any metric space $X$ one can study quasi-isometric embeddings of $\R^{n}$ 
into $X$, these are called \emph{quasiflats of rank $n$}.\index{quasiflats}
Generalizing the above Euclidean example, one way to quantify 
non-hyperbolicity of a metric space $X$ is by calculating its \emph{geometric 
rank}\index{geometric rank}, namely the largest 
rank of any quasiflat in $X$. Having geometric rank greater than one is 
an obstruction to hyperbolicity, as it forces $X$ to have arbitrarily 
fat triangles.

We now record a technical lemma which will be useful in 
\fullref{sectionhyperbolicitylowcomplexity}.
The proof is an exercise in manipulating thin triangles (see \cite{Behrstock:thesis} 
for a proof).  When the ambient space is clear we use the notation 
$B_{\delta}(x)$ to denote the closed ball of radius $\delta$ around a 
point $x$.

\begin{lem}\label{lemmaeasythintriangleargument}
    Let $Z$ be a $\delta$--hyperbolic space and $\gamma$ a geodesic in 
    $Z$. If $\phi\co Z\to 2^{\gamma}$ sends each point of $Z$ to the set of 
    closest points on $\gamma$, then for each point $z\in Z$ we have 
    $\diam(\phi(z))\leq 4\delta+2$.
\end{lem}

\subsection{Asymptotic Cones}

Ultrafilters and ultraproducts have been used in logic since G\"{o}del
introduced them in the 1930's. Among other applications they can be
used to give a slick proof of the compactness theorem of first order
logic and to rigorize calculus, {\it a la} Robinson's non-standard
analysis. Van~den~Dries and Wilkie introduced these tools to
topologists when they reformulated (and slightly strengthened)
Gromov's famous Polynomial growth Theorem using this language (compare
Gromov \cite{Gromov:PolynomialGrowth} and van den Dries and Wilkie
\cite{DriesWilkie}). These techniques streamlined parts of the proof,
replacing an iterated procedure of passing to subsequences with one
ultralimit, and broadened the applicability of Gromov's construction,
Gromov's original definition of asymptotic cone only applied to
nilpotent groups whereas the ultrafied definition works for any
finitely generated group. (For other applications of asymptotic cones
in geometric group theory see: Burillo \cite{Burillo:rankBS}, Bridson
\cite{Bridson:Cones}, Dru{\c{t}}u\cite{Drutu:QIinvariantsCones},
Dru{\c{t}}u and Sapir \cite{DrutuSapir:TreeGraded}, Kapovich and Leeb
\cite{KapovichLeeb:haken} Kleiner and Leeb
\cite{KleinerLeeb:buildings}, Kramer, Shelah, Tent and Thomas
\cite{KSTT} and Riley \cite{Riley:ConesConnectedness}.)

 A \emph{non-principal ultrafilter} on the integers, 
 denoted $\omega$, is a nonempty collection of 
 sets of integers with the following 
 properties:
 \begin{enumerate} 
     \item If $S_{1}\in\omega$ and $S_{2}\in\omega$, then 
     $S_{1}\cap S_{2}\in\omega$.
     \item If $S_{1}\subset S_{2}$ and $S_{1}\in\omega$, then 
     $S_{2}\in\omega$.
     
     \item It is \emph{maximal} in the sense that for each  
      $S\subset \Z$ exactly one of the following must occur:
      $S\in\omega$ or $\Z\setminus S\in\omega$.
      
      \item $\omega$ does not contain any finite sets. (This is the 
      \emph{non-principal} aspect.)
 \end{enumerate}

 Since we are only concerned with non-principal ultrafilter, we will 
 simply refer to these as \emph{ultrafilters}.

 An equivalent way to view an ultrafilter is as a finitely additive 
 probability measure, with no atoms, 
 defined on all subsets of integers and which takes 
 values in $\{0,1\}$. This way of thinking is 
 consistent with the intuition that the ultrafilter sees the 
 behavior on only the large subsets of $\N$ and when this is done one 
  says that $\omega$--almost every integer has a given 
 property when the set of integers with this property is in $\omega$.

 For an ultrafilter $\omega$, a metric space $(X, d)$, and a sequence of 
 points $\<y_{i}\>_{i\in\N}$ we define $y$  to be 
 the \emph{ultralimit of $\<y_{i}\>_{i\in\N}$ with respect to 
 $\omega$}, denoted $y=\ulimi y_{i}$, 
 \index{ultralimit}\index{$\ulim$} 
 if and only if for all $\epsilon>0$ one has 
 $\{i\in\N:d(y_{i},y)<\epsilon\}\in\omega$.
 When  the ultralimit of a sequence exists it is a 
 unique point; further, 
 it is easy to see that when the sequence $\<y_{i}\>_{i\in\N}$ 
 converges, then the ultralimit is the same as the ordinary limit.
 It is also worth noting that in a compact metric space 
 ultralimits always exist.

 Fix an ultrafilter $\omega$ and a family of based metric spaces 
 $(X_{i}, x_{i}, d_{i})$. 
 Using the ultrafilter, a pseudo-distance on $\prod_{i\in\N}X_{i}$ is 
 provided by: 
 $$d_{\omega}(\<a_{i}\>,\<b_{i}\>)=\ulimi 
  d_{X_{i}}(a_{i},b_{i})\in [0,\infty].$$
 \noindent
 The \emph{ultralimit of $(X_{i},x_{i})$}\index{$\ulim$} 
 is then defined to be:
  $$\ulimi (X_{i},x_{i}) =\{y\in\prod_{i\in\N}X_{i}: 
  d_{\omega}(y,\<x_{i}\>)<\infty\}/ \sim,$$
  where for two points  $y,z\in \prod_{i\in\N}X_{i}$ we define $y\sim z$ 
  if and only if  $d_{\omega}(y,z)=0$.\index{$\sim$}
  
  The pseudo-metric $d_{\omega}$ takes values in 
  $[0,\infty]$---considering only those points  in\break $\ulim (X_{i},x_{i})$ 
  which are a 
  finite distance from the basepoint restricts our attention to a 
  maximal subset where the pseudo-distance does not obtain 
  the value $\infty$. The relation $\sim$ quotients out points whose 
  pseudo-distance from each other is zero;  in an analogy to measure theory 
  we think of this process as identifying sequences  
  which agree almost everywhere.
  These two conditions combine to make $d_{\omega}(y,z)$ into 
  a metric.
 
  Just as ultralimits for sequences of 
  points generalize the topological notion of limit, ultralimits of 
  sequences of metric spaces generalize Hausdorff limits of metric 
  spaces. In particular, for a Hausdorff precompact family 
  of metric spaces the ultralimit of the sequence is a 
  limit point with respect to the Hausdorff topology (see 
  \cite{KleinerLeeb:buildings} for a proof of this and other related facts).

 \begin{defn} The \emph{asymptotic cone of $(X,\<x_{i}\>,\<d_{i}\>)$ 
     relative to the 
 ultrafilter $\omega$}\index{asymptotic cone} is defined by:
 $$\coneomega (X,\<x_{i}\>,\<d_{i}\>)=\ulimi 
 (X,x_{i},\frac{1}{d_{i}}).$$
  When it is not a source of confusion we tend to suppress writing the 
 basepoint.
 Note that for for homogeneous spaces (like Cayley graphs), the
 asymptotic cone does not depend on the basepoint $\<x_{i}\>$ and thus 
 we shall always take this to be $\<1\>$ and drop it from the notation. 
\end{defn}

When taking such ultralimits, we often use the notation $\ulimif 
(x_{\alpha}, 0_{\alpha})$, when there is an implied ambient space 
$T_{\alpha}$ with basepoint $0_{\alpha}$, in order to emphasize that we 
are looking at an  ultralimit in $\cone (T_{\alpha}, 0_{\alpha})$.

When considering the asymptotic cone, especially for the first time, 
the following quotation from Thomas Hobbes (1588--1679), 
will likely seem relevant.
\begin{quotation}
    To understand this for sense it is not required that a man should be 
    a geometrician or a logician, but that he should be mad.
\end{quotation}

Asymptotic cones provide a way to replace a metric space with 
a ``limiting'' space which carries information about sequences in 
the original space which leave every compact set. This process 
encodes the asymptotic geometry of a space  into  
standard algebraic topology invariants of its asymptotic cone.

A common use of asymptotic cones is in their relation to 
hyperbolicity as 
demonstrated by the next result and \fullref{treeimplieshyper}. 
First we give a preliminary definition.

\begin{defn} An $\R$--tree\index{$\R$--tree} is a metric space $(X,d)$ such that between 
any two points $a,b\in X$ there exists a unique topological arc $\gamma$
connecting them and $\gamma$ is isometric to the interval 
$[0,d_{X}(a,b)]\subset \R$.
\end{defn}

\begin{prop}\label{hyperimpliestree} For a sequence $\delta_{i}\to 0$, 
    the ultralimit of $\delta_{i}$--hyperbolic spaces is an $\R$--tree.
    
    In particular, if $X$ is a $\delta$--hyperbolic 
space, then $\cone X$ is an $\R$--tree.
\end{prop}

If $X$ is a $\delta$--hyperbolic space which is sufficiently
complicated (eg, the cayley graph of a non-elementary hyperbolic
group) then  $\cone X$ is an $\R$--trees with uncountable branching at 
every point.

A partial converse to \fullref{hyperimpliestree} is provided by 
the following well known result. The statement of this 
theorem first appeared, without proof, in Gromov \cite{Gromov:hyperbolic}.  
A proof can be found in Dru\c{t}u
\cite{Drutu:QIinvariantsCones}.

\begin{thm}\label{treeimplieshyper}
    If every asymptotic cone of a metric space $X$ is an $\R$--tree,
    then $X$ is $\delta$--hyperbolic.  
\end{thm}

The following is an easy observation that simplifies the situation 
when  
dealing with finitely many equivalence classes of $\seq$, where $\seq$ 
is the set of sequences in a countable alphabet considered up to the
equivalence relation they they agree $\omega$--almost everywhere. When 
we consider $\seq$ in the sequel the alphabet will consist of
subsurfaces of a given surface.

\begin{lem}\label{easyultralimitobservation}
    For any finite set $\Gamma\subset\seq$,  
    the elements of $\Gamma$ are pairwise distinct in $\seq$ if and 
    only if there  exists a set $K\in\omega$ 
    for which each $\gamma,\gamma'\in\Gamma$ has 
    $\gamma_{i}\neq\gamma_{i}'$ for every $i\in K$
\end{lem}

\begin{proof}  Fix a finite set $\Gamma\subset\seq$.

    {\bf ($\Longrightarrow$)}\qua  The maximal clause in the definition of
    ultrafilter states that for each $K\subset\N$ either
    $K\in\omega$ or $\N\setminus K\in\omega$.
    We suppose that for
    each pair of elements $\gamma,\gamma'\in\Gamma$
    we have $\gamma\nsim\gamma'$ and thus the set $K$ of indices where
    $\gamma_{i}=\gamma_{i}'$  must have
    $K\notin\omega$. Maximality of the ultrafilter then implies
    that $\N\setminus K\in\omega$.
    Thus  for each $i\in\N\setminus K\in\omega$
    we  have  $\gamma_{i}\neq\gamma_{i}'$, we define
    $\N\setminus K=K'_{\gamma,\gamma'}$.

    Since ultrafilters are closed under
    finite intersections, the intersection over all pairs $\gamma,\gamma'$
    of $K_{\gamma,\gamma'}$ yields a set $J\in\omega$ where
    $\gamma_{i}\neq\gamma_{i}'$ for every $i\in J$ and every
    $\gamma,\gamma'\in\Gamma$.

    {\bf ($\Longleftarrow$)}\qua
    Let  $K\in\omega$ be the set of indices for which $\gamma_{i}\neq\gamma_{i}'$
    for each $i\in K$ and $\gamma,\gamma'\in\Gamma$.
 Then for any $\gamma, \gamma'$ we have
   $\{i: \gamma_{i} = \gamma'_{i}\}\in\N\setminus K$ and hence is not
   in $\omega$. Thus $\gamma\nsim\gamma'$.
\end{proof}

We end this section with a summary of some standard results about
asymptotic cones that we will use in the sequel (see 
Bridson and Haefliger \cite{BridsonHaefliger}, Kleiner and Leeb
\cite{KleinerLeeb:buildings} or Kapovich \cite{Kapovich:GGTbook}).

\begin{prop}\label{generalconefacts} Fix a 
    non-principal ultrafilter $\omega$.
\begin{enumerate}
    
    \item $\con(X)$ is a complete metric space.
    
    \item $\con(X_{1}\times X_{2})=\con(X_{1})\times\con(X_{2})$

    \item $\con\R^{n}=\R^{n}$

    \item The asymptotic cone of a geodesic space is a geodesic space.
    
    \item A $(K,C)$--quasi-isometry between metric spaces induces a 
    $K$--bi-Lipschitz map between their asymptotic cones.
\end{enumerate}

\end{prop}

\section{Background II: Surfaces and related structures}
\label{chaptersurfacesbackground}
\subsection{Complex of Curves}\label{ccbackground}

We use $S$ to denote a connected, orientable surface of genus
$g=g(S)$ with $p=p(S)$ punctures. (Note that we will interchangeably 
use the 
words puncture and boundary, as the distinction does not affect the 
results contained in this work.) We use the terms \emph{subsurface} and 
\emph{domain}
to refer to a homotopy class of an essential, 
non-peripheral, connected subsurface of $S$ (subsurfaces are not 
assumed to be proper unless explicitly stated). 
When we 
refer to the boundary of a domain this will mean the collection of 
(homotopy classes of) 
non-peripheral closed curves which bound the domain as a subset of $S$.

We now recall the construction of the complex of curves and
relevant machinery developed by Masur and Minsky which we will use in 
our study; for further details consult \cite{MasurMinsky:complex2}. 
We use $\xi(S)=3g(S) +p(S)-3$\index{$\xi(S)$} to 
quantify the complexity\footnote{In 
\cite{MasurMinsky:complex2} they use $\xi(S)=3g(S)+p(S)$, 
but we use the current version because it has better additive 
properties when generalized to disconnected 
surfaces, as we need in \fullref{sectiontransversaltree}.} of the surface $S$.
Recall that when positive, $3g(S)+p(S)-3$  
is the maximal number of disjoint homotopy classes of 
essential and non-peripheral   
simple  closed curves which can be simultaneously realized on $S$.
Naturality of $\xi(S)$ as a measure of complexity is justified by the
property that it decreases when one passes from a surface to a 
proper subsurface (recall our convention of considering surfaces up to 
homotopy). Since our interest is in hyperbolic surfaces and their 
subsurfaces, we only 
consider subsurfaces with $\xi>-2$ (thus ignoring the disk and the 
sphere). Additionally, as it is not a hyperbolic surface nor does it 
appear as a subsurface of any hyperbolic surface we will usually ignore 
$S_{1,0}$ (although much of our discussion has analogues for this case); 
thus $\xi=0$ is used only to denote the thrice punctured sphere.

Introduced by Harvey \cite{Harvey:Boundary} to study the boundary of 
Teichm\"{u}ller space, the complex of curves has proven to be a 
useful tool in the study of Teichm\"{u}ller space, mapping class 
groups, and $3$--manifolds. The complex of curves is a finite 
dimensional complex which encodes information about the surface via the 
combinatorics of simple closed curves. Analysis using the curve 
complex is necessarily delicate since the complex is 
locally infinite except in a few low genus cases. 

\begin{defn}
The \emph{complex of curves for $S$}\index{curve complex}, 
denoted $\C(S)$, \index{$\C(S)$}  
consists of a \emph{vertex} for 
every  homotopy class of a simple closed curve which is both 
non-trivial and non-peripheral. The \emph{$N$--simplices} of $\C(S)$ are given by 
collections of $N+1$ vertices whose homotopy classes can all 
simultaneously be realized disjointly on $S$.
\end{defn}

This definition works well when $\xi(S)>1$, however it must be modified slightly 
for the surfaces of small complexity, which we refer to as the 
\emph{sporadic cases}\index{curve complex!sporadic cases}. 
The two cases where $\xi(S)=1$ 
are $S_{1,1}$ and 
$S_{0,4}$ (the following discussion works for $S_{1,0}$ as well); 
here any two distinct homotopy classes of essential simple closed 
curves must intersect (at least once in  $S_{1,1}$ and twice in 
$S_{0,4}$). In these cases the above definition of $\C(S)$ would have 
no edges, only vertices. Accordingly, we modify the definition so that 
an edge is added between two homotopy classes when they intersect the 
minimal amount possible on the surface (ie, once for  $S_{1,1}$ and 
twice for  $S_{0,4}$). With this definition 
$\C(S_{1,1})$ and $\C(S_{0,4})$ are each connected, indeed 
they are each isometric to the classical Farey graph.

When $\xi(S)=0, -2,$ or $-3$ then $\C(S)$ is empty. The final modification we 
make is for the case where we have $A\subset S$ with $\xi(A)=-1$, the annulus. 
The annulus doesn't support a finite area hyperbolic metric, 
so our interest in it derives from the fact that it arises
as a subsurface of hyperbolic surfaces. Indeed, annuli will play a 
crucial role as they will be used to capture information about Dehn 
twists. 
Given an annulus $A\subset S$, we define $\C(A)$ to be based homotopy classes of 
arcs connecting one boundary component of the annulus to the other. 
More precisely,  
denoting by $\widetilde{A}$ the annular cover of $S$ to which $A$ lifts 
homeomorphically, we use the compactification of $\H^{2}$ as the 
closed unit disk to obtain a closed annulus $\widehat{A}$. We define 
the vertices of $\C(A)$ to be homotopy classes of paths  connecting one 
boundary component of  $\widehat{A}$ to the other, where the homotopies 
are required to fix the endpoints. Edges of $\C(A)$ are pairs of 
vertices which have representatives with disjoint interiors. Giving 
edges a Euclidean metric of length one, in \cite{MasurMinsky:complex2} it 
is proven that $\C(A)$ is quasi-isometric to $\Z$.

The following foundational result, which is the main theorem of
\cite{MasurMinsky:complex1}, will be useful in our later analysis. For
another proof see Bowditch \cite{Bowditch:Curvecomplex}, which provides a
constructive computation of a bound on the hyperbolicity constant.

\begin{thm}[Hyperbolicity of $\C(S)$ \cite{MasurMinsky:complex1,Bowditch:Curvecomplex}] \label{curvecomplexhyperbolic}
    For any surface $S$, $\C(S)$ is an infinite diameter
    $\delta$--hyperbolic space (as long as it is non-empty). 
\end{thm}

Throughout this paper we use the convention that ``intersection'' 
refers to transverse intersection. Thus for example if we consider a 
subsurface $Y\subsetneq S$ and an element $\gamma\in\partial Y$ which 
is not homotopic to a puncture of $S$ then we consider the annulus around 
$\gamma$ to not intersect $Y$.

For a surface with punctures one can consider the \emph{arc
complex}\index{arc complex} $\C'(S)$, \index{$\C'(S)$} which is a
close relative of the curve complex.  When $\xi(S)>-1$ we define the
vertices\footnote{This definition differs from that given by Harer
\cite{Harer:vcd} where an arc complex is considered consisting of only
arcs and thus does not contain $\C(S)$ as a subcomplex. Our definition
agrees with that in \cite{MasurMinsky:complex2}.}  of $\C'(S)$
(denoted $\C'_{0}(S)$\index{$\C'_{0}(S)$}) to consist of elements of
$\C_{0}(S)$ as well as homotopy classes of simple arcs on $S$ with
endpoints lying on punctures of $S$, which don't bound a disk or a
once punctured disk on either side.  As done for the curve complex, we
define $N$--simplices of $\C'(S)$ to be collections of $N+1$ vertices
which can simultaneously be realized on the surface as disjoint arcs
and curves.  For annuli we define $\C'(A)=\C(A)$.

The arc complex arises naturally when one tries to ``project'' an element 
$\gamma\in\C(S)$ into $\C(Y)$ where $Y\subset S$. When $\xi(Y)>0$ we define 
$$\pi'_{Y}\co\C_{0}(S)\to 2^{C'_{0}(Y)}$$ 
by the following:
\begin{itemize}
    \item If $\gamma\cap Y=\emptyset$, then define $\pi'_{Y}(\gamma)=\emptyset$.

    \item If $\gamma\subset Y,$ then define $\pi'_{Y}(\gamma)=\{\gamma\}$.

    \item If $\gamma\cap\partial Y\neq \emptyset$, then after putting 
    $\gamma$ in a position so it has minimal intersection with 
    $\partial Y$ we identify
    parallel arcs of $\gamma\cap Y$ and define  $\pi'_{Y}(\gamma)$ 
    to be the union of these arcs and any closed curves in $\gamma\cap\partial 
    Y\neq \emptyset$.
\end{itemize}
In the last case, it follows 
from the definition of $\C'(Y)$ that $\pi'_{Y}(\gamma)$ is a subset of 
$\C'_{0}(Y)$ with diameter at most one. Thus whenever 
$\pi'_{Y}(\gamma)\neq\emptyset$, it  
is a subset of $\C'(Y)$ of diameter at most one. 
Moreover, when $\xi(S)>0$, 
there is a map $\phi_{Y}\co\C'_{0}(Y)\to 2^{\C_{0}(Y)}$, which implies 
that the arc complex embeds into the curve complex as a cobounded 
subset---the map $\phi_{Y}$ sends each arc to the 
boundary curves of a regular neighborhood of its union with $\partial Y$. 
For $Y\subset S$ with $\xi(Y)>0$ we define
$\pi_{Y}=\phi_{Y}\circ\pi_{Y}'\co  
\C(S)\to \C(Y)$\index{$\pi_{X}$} (this map is actually to the power 
set of $\C(Y)$, but for the sake of readability and since the map is to 
bounded sets, we abuse notation and simply write  $\C(Y)$).

When $\xi(Y)=-1$ any curve $\gamma$ which crosses $Y$ 
transversally has a lift $\widetilde{\gamma}\in\widetilde{Y}$ with at 
least one component which connects the two boundary components of 
$\widehat{Y}$. Together, the collection of lifts which connect the 
boundary components form a finite subset of $\C'(Y)$ with diameter at most 
1, define $\pi_{Y}(\gamma)$ to be this set. If $\gamma$ doesn't 
intersect $Y$ or is the core curve of $Y,$ then define 
$\pi_{Y}(\gamma)=\emptyset$. Also, for consistency  we 
define $\pi_{Y}\co v\to\{v\}$ for $v\in\C_{0}(Y)$.

Since we often work with subsets of $\C(Y)$ rather than points, we 
use the following notation. \index{$d_{X}$} \index{$\diam_{X}$}
For a set valued map $f\co X\to 2^{Y}$ and a set $A\subset X$ we 
define $f(A)=\cup_{x\in A} f(x)$.
If $X\subseteq S$ and $\mu,\nu\in\C(S)$ we adopt the notation    
$$d_{X}(\mu,\nu)=d_{\C(X)}(\pi_{X}(\mu),\pi_{X}(\gamma)).$$
Further, given sets $A,B\in\C(S)$ we set 
$d_{X}(A,B)=\min\{d_{X}(\alpha,\beta):\alpha\in A \mbox{ and } 
\beta\in B\}$. Also, we write $\diam_{X}(A)$ to refer to the diameter of the set
$\pi_{\C(X)}(A)$ and $\diam_{X}(A,B)$ for
$\diam_{\C(X)}(A\cup B)$ in order to emphasize the symmetry between our use of
minimal distance and diameter. (These same conventions apply when 
considering markings, as defined in the next section.)

An extremely useful result concerning these projections is the 
following (see \cite{MasurMinsky:complex2} for the original proof, and
\cite{Minsky:ELC1} where the bound is corrected from 2 to 3):

\begin {lem}[Lipschitz projection  
\cite{MasurMinsky:complex2}]\label{MMtwlpl}
Let $Y$ be a subdomain of $S$.  For any simplex $\rho$ in $\C(S)$, if
$\pi_{Y}(\rho) \neq \emptyset$ then $\diam_{Y}(\rho) \leq 3$.
\end{lem}

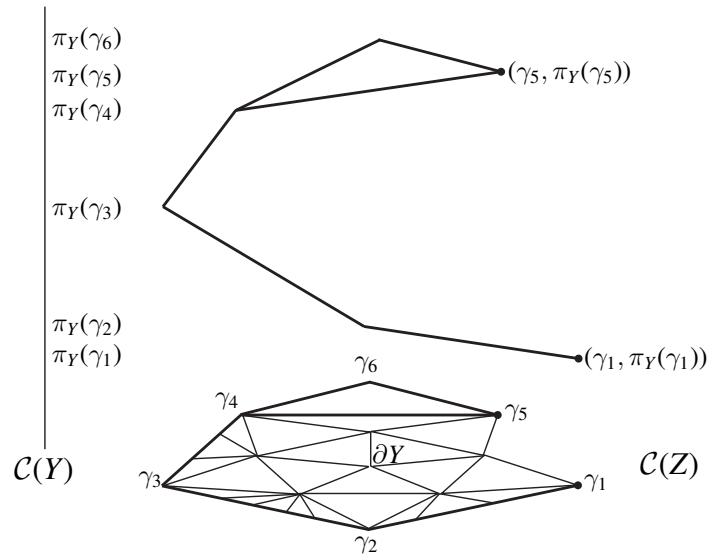
\begin{figure}[ht!]
\labellist\small
\pinlabel* {$\gamma_{3}$} [r] at 55 27
\pinlabel {$\pi_{Y}(\gamma_{5})$} [l] at 13 170
\pinlabel* {$\gamma_{4}$} [br] at 82 52
\pinlabel {$\pi_{Y}(\gamma_{4})$} [l] at 13 158
\pinlabel {$\gamma_{5}$} [l] at 173 51
\pinlabel {$\pi_{Y}(\gamma_{3})$} [l] at 13 123
\pinlabel {$\gamma_{1}$} [l] at 201 26
\pinlabel {$\pi_{Y}(\gamma_{1})$} [l] at 13 71
\pinlabel {$\gamma_{2}$} [t] at 127 10
\pinlabel {$\pi_{Y}(\gamma_{2})$} [l] at 13 82
\pinlabel {$\gamma_{6}$} [b] at 127 62
\pinlabel {$\pi_{Y}(\gamma_{6})$} [l] at 13 183
\pinlabel {$(\gamma_{1},\pi_{Y}(\gamma_{1}))$} [l] at 201 70
\pinlabel {$(\gamma_{5},\pi_{Y}(\gamma_{5}))$} [l] at 174 171
\pinlabel {$\partial Y$} [bl] <0pt,-2pt> at 127 33
\pinlabel {\large{$\C(Y)$}} [t] at 13 38
\pinlabel {\large{$\C(Z)$}} at 234 32
\endlabellist
\centerline{\includegraphics[width=3.5truein]{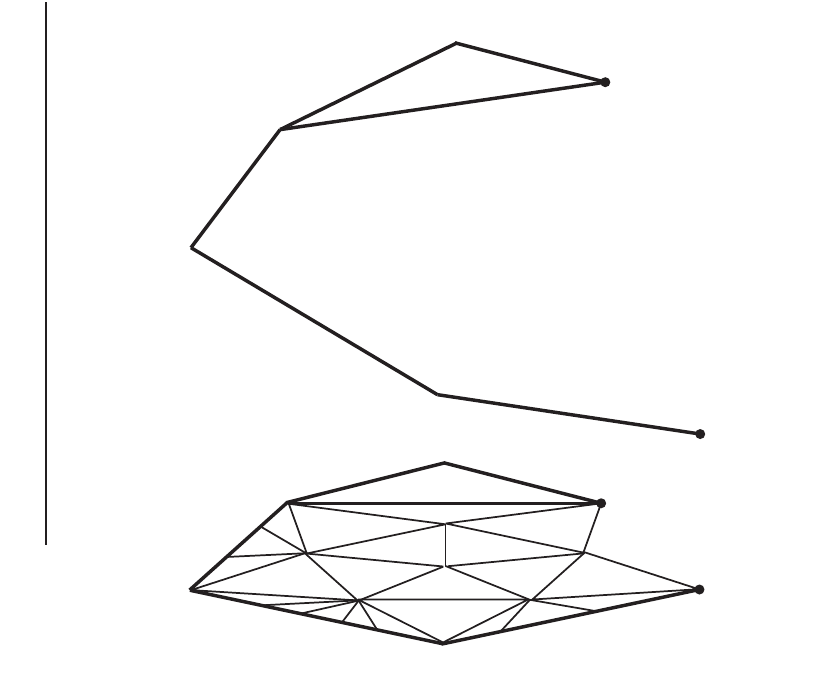}}
\caption[Graph of $\pi_{Z\to Y}$]{Letting  
 $Y\subset Z$, the above is a caricature of the graph of 
 $\pi_{Z\to Y}\co \C(Z)\setminus B_{1}(\partial Y)\to \C(Y)$  near the 
 point $\partial Y\in \C(Z)$.}
\label{nestedprojpicts}
\end{figure}

In light of this lemma, for any pair of domains 
$Y\subset Z$ we consider  
 $\C(Z)\setminus B_{1}(\partial Y)$ with the path metric, ie, the 
 distance between $\gamma,\gamma'\in\C(Z)\setminus B_{1}(\partial Y)$ 
 is the length of the shortest path connecting them in $\C(Z)\setminus 
 B_{1}(\partial Y)$. With this metric there is a
 \emph{coarsely Lipschitz}\index{coarsely Lipschitz} map 
 $\pi_{Z\to Y}\co\C(Z)\setminus B_{1}(\partial Y) \to 
 \C(Y)$\index{$\pi_{Z\to Y}$}, where a 
 map is defined to be 
  \emph{$(K,C)$--coarsely 
Lipschitz} if and only if there exist a pair of constants $K,C$ 
such that for each $a,b\in \C(Z)\setminus B_{1}(\partial Y)$ we have 
$d_{Y}(\pi_{Z\to Y}(a),\pi_{Z\to Y}(b))\leq K d_{\C(Z)\setminus B_{1}(\partial 
Y)}(a,b) +C$. By the above lemma and the fact that 
$\pi_{Z\to Y}$ sends points of $\C(Z)\setminus B_{1}(\partial Y)$ 
to subsets of $\C(Y)$ of diameter at most 3, 
we have the following which we shall use in 
 \fullref{chapterProjectionEstimates}. (See \fullref{nestedprojpicts} 
 for a cartoon of the behavior of this map near $\partial Y$.)

\begin{cor}\label{projectioncoarselylip} 
For any domain $Y\subset Z$, endowing $\C(Z)\setminus B_{1}(\partial Y)$  
with the path metric, we have  
$$\pi_{Z\to Y}\co\C_{1}(Z)\setminus B_{1}(\partial Y)\to \C(Y)$$
is coarsely Lipschitz (with constants $K=3$ and $C=3$).
\end{cor}

In \fullref{nestedprojpicts}, we observe that although $\pi_{Z\to Y}$ 
is Lipschitz, pairs of points like $\gamma_{1}$ and $\gamma_{5}$, 
although distance 2 in $\C(Z)$, can be found which are arbitrarily far apart 
in the path 
metric on $\C(Z)\setminus B_{1}(\partial Y)$ and thus their distance 
when projected to $\C(Y)$ are also made arbitrarily large.

Letting  $\C_{1}(S)$ denote the one-skeleton of  $\C(S)$, we state 
the following theorem which 
is one of the key technical results from \cite{MasurMinsky:complex2}. 
\begin{thm}[Bounded geodesic image \cite{MasurMinsky:complex2}]\label{MMtwbgi}
    Let $Y\subsetneq S$ with $\xi(Y)\neq 0$ and let $g$ be a geodesic 
    segment, ray, or line in $\C_{1}(S)$, such that $\pi_{Y}(v) 
    \neq\emptyset$ for every vertex $v$ of $g$.
    There is a constant $D$ depending only on $\xi(S)$ so that 
    $\diam_{Y}(g)\leq D$.
\end{thm}

\subsection{Markings}\label{markingbackground}

In this section we describe the quasi-isometric model we will use for
the mapping class group and explain the tools developed in 
\cite{MasurMinsky:complex2} for computing with this model.

\begin{figure}[ht!]\label{figurecompletecleanmarking}
    \labellist\small
    \pinlabel {$\gamma_{1}$} [t] at 52 11
    \pinlabel {$\gamma_{2}$} [t] <0pt, -1pt> at 125 46
    \pinlabel* {$t_{1}$} [tr] at 22 42
    \pinlabel {$t_{2}$} [b] at 176 70
    \endlabellist    
\centerline{\includegraphics[width=2.5truein]{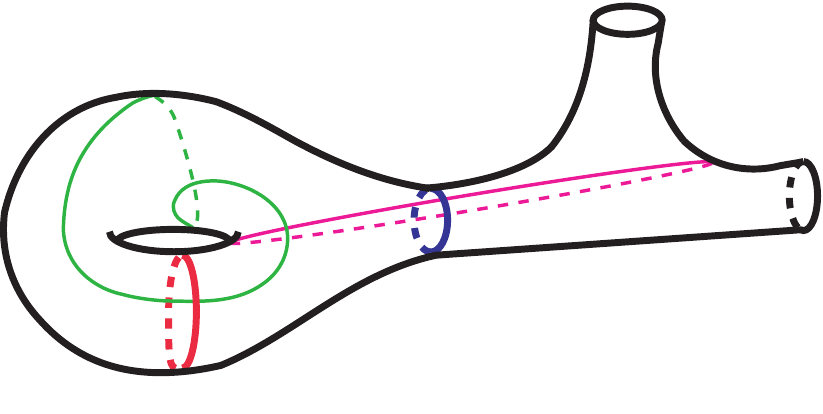}}
\caption[A complete clean marking]{A complete clean marking 
$\mu\in\M(S_{1,2})$ with $\base(\mu)=\{\gamma_{1}, \gamma_{2}\}$}
\end{figure}
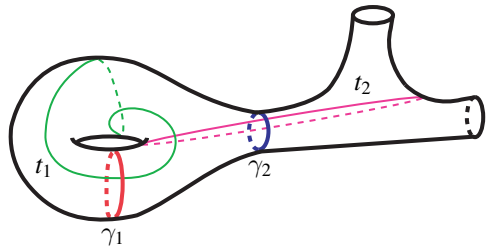

\begin{defn} A \emph{marking}\index{marking}, $\mu$, on $S$ is a collection of 
    \emph{base curves} to each of which we (may) associate a 
    \emph{transverse curve}. These collections are made subject to the 
    constraints:
    \begin{itemize}
        \item The \emph{base curves}\index{marking!base curves},   
        $\base(\mu)=\{\gamma_{1},\ldots,\gamma_{n}\}$, form a simplex 
        in $\C(S)$. \index{$\base(\mu)$}
    
        \item 
	The \emph{transverse curve}\index{marking!transverse curve}, 
	$t$, associated to a given base 
        curve $\gamma$ is either empty or an element of $\C_{0}(S)$ which 
        intersects $t$ once (or twice if $S=S_{0,4}$)  
	and projects to 
	a subset of diameter at most 1 in the annular complex $\C(\gamma)$.
    \end{itemize}
    
    When the  transverse curve $t$ is empty  
     we say that $\gamma$ doesn't have a transverse curve. 
    
    If 	the simplex formed by $\base(\mu)$ in $\C(S)$ is top 
    dimensional and every curve has a transverse curve,
    then we say the marking is \emph{complete}.\index{marking!complete}
    
    If a marking has the following two properties, then
    we say the marking is \emph{clean}\index{marking!clean}. First, for
    each $\gamma\in\base{\mu}$ its transversal $t$ 
    is disjoint from the rest of the base curves. Second, for each 
    $\gamma$ and $t$ as above,  $t\cup\gamma$ fills a
    surface denoted $F(t,\gamma)$ with $\xi(F(t,\gamma))=1$ and in which 
    $d_{\C(F(t,\gamma))}(t,\gamma)=1$. (See \fullref{figurecompletecleanmarking}.)
\end{defn}

Let $\mu$ denote a complete clean marking with pairs $(\gamma_{i}, 
t(\gamma_{i}))$, 
we take as \emph{elementary moves}\index{marking!elementary moves} 
the following two relations on the set 
of complete clean markings:
\begin{enumerate}
    \item \emph{Twist}: For some $i$, we replace $(\gamma_{i}, t(\gamma_{i}))$ by 
    $(\gamma_{i}, t'(\gamma_{i}))$ where $t'(\gamma_{i})$ is the result 
    of one full (or half when possible) twist of $t(\gamma_{i})$
    around $\gamma_{i}$. The rest of the pairs are left unchanged.

    \item \emph{Flip}: For some $i$ we swap the roles of the base and transverse 
    curves between $\gamma_{i}$ and $t(\gamma_{i})$. After 
    doing this the complete marking may no longer be clean, one then 
    needs to replace the new marking with a compatible clean one.
\end{enumerate}

We say the clean marking $\mu$ is \emph{compatible} with $\mu'$ 
if they have the same base 
curves, each base curve $\gamma$ has a transverse curve $t$ in one 
marking if and only if 
it has a transverse curve $t'$ in the other marking, and when  
$t$ exists then  $d_{\gamma}(t,t')$ is minimal among all choices 
of $t'$.

In \cite{MasurMinsky:complex2} it is shown that there exists a bound 
(depending only on the topological type of $S$) on the number of clean 
markings which are compatible with any other given marking. Thus even 
though the flip move is defined by choosing an arbitrary compatible 
complete clean marking, it is canonical up to some uniformly bounded amount 
of ambiguity.

One then defines the \emph{marking complex}\index{marking!complex of}, 
denoted $\M(S)$\index{$\M(S)$}, to be 
the graph formed by taking the complete clean markings on $S$ as 
vertices and connecting two vertices by an edge if they differ by an 
elementary move. It is not hard to check that $\M(S)$ is a locally 
finite graph and that the mapping class group acts on it 
cocompactly and 
properly discontinuously. A consequence of this which we use throughout 
this work is:

\begin{cor}{\rm\cite{MasurMinsky:complex2}}\qua $\M(S)$ is quasi-isometric to 
    the mapping class group.\index{mapping class group}
\end{cor}

For any subsurface $Y\subset S$ with $\C(Y)\neq \emptyset$ 
we considered in the previous section the \emph{subsurface projection}, 
$\pi_{Y}\co\C(S)\to 2^{\C(Y)}$.\index{$\pi_{X}$}
More generally one can consider subsurface projections from the  
marking complex; since this definition generalizes the above map we also
denote it $\pi_{Y}\co\M(S)\to 2^{\C(Y)}$. When $\mu\in\M(S)$ we define 
$\pi_{Y}(\mu)$ to be $\pi_{Y}(\base(\mu))$ unless $Y$ is an annulus 
around a curve $\gamma\in\base(\mu)$, in which case we 
define $\pi_{Y}(\mu)=\pi_{Y}(t)$ where $t$ is the transverse curve to 
$\gamma$. An 
important observation is that when $\mu$ is a complete marking on $S$ 
then $\pi_{Y}(\mu)\neq\emptyset$ for each $Y\subset S$.

There is also a Lipschitz projection property for projections of 
markings.

\begin {lem}[Elementary move projection \cite{MasurMinsky:complex2}]
\label{MMtwemp} 
If $\mu,\mu'\in\M(S)$ differ by one elementary move, then for any 
domain $Y\subset S$ with $\xi(Y)\neq 0$, $$d_{Y}(\mu,\mu') \leq 4.$$
\end{lem}

\subsection{Hierarchies}\label{hierarchiesbackground}

Thurston's classification theorem for surface homeomorphisms 
\cite{Thurston:surfaces} gives a layered structure for homeomorphisms 
based upon studying subsurfaces which are preserved (perhaps under an 
iterate of the homeomorphism). In \cite{MasurMinsky:complex2}, Masur and Minsky
refine and further elucidate  a layered structure for the mapping 
class group. Using the marking complex and the complex of curves, they 
provide a way to compare how relatively complicated two given mapping 
class group elements are on any subsurface (not just ones that are 
eventually periodic!). Towards this end, they introduce an object called 
a {\it hierarchy} which ties all the subsurface comparisons together. 

\begin{defn} Let $Y\subset S$  with $\xi(Y)>1$.  A sequence of simplices 
$v_{0}, v_{1},\ldots ,v_{n}$ is called 
\emph{tight}\index{tight!sequence} if:
\begin{enumerate}
    \item For each $0\leq i<j\leq n$, and vertices
    $v_{i}'\subset v_{i}$ and  $v_{j}'\subset v_{j}$,   
    $$d_{Y}(v_{i}',v_{j}')=|i-j|.$$

    \item For each $0<i<n$ we have that $v_{i}$ is the  
    boundary  of the subsurface filled by $v_{i-1}\cup v_{i+1}$.
\end{enumerate}

When $\xi(Y)=1$ we consider a sequence to be \emph{tight} if and only 
if it is the vertex sequence of a geodesic.

When $\xi(Y)=-1$ we consider a sequence to be \emph{tight} if and only 
if it is the vertex sequence of a geodesic where the set of endpoints 
on $\partial \widehat{Y}$ of arcs representing the vertices equals the set 
of endpoints of the first and last arc.
\end{defn}

We often use the following decorated version of a tight sequence.
\begin{defn}  A  \emph{tight geodesic}\index{tight!geodesic} $g$ in $\C(Y)$ consists of 
    a tight sequence 
    $v_{0}, v_{1}\ldots, v_{n}$ and a pair of markings $\I=\I(g)$ and 
    $\T=\T(g)$ (called the \emph{initial} and \emph{terminal markings} 
    for $g$)\index{marking!initial}\index{marking!terminal} 
    such that $v_{0}$ is a subset of the simplex $\base(\I)$ and $v_{n}$ 
    is a subset of  $\base(\T)$.
    
    The integer $n$ is called the length of $g$. The \emph{domain} (sometimes 
    called the \emph{support}) of $g$ refers to the surface $Y,$ written 
    $D(g)=Y$.

\end{defn}

Below we explain a relationship between tight geodesics occurring in 
different subsurfaces.

When $\mu\in\M(S)$ we write 
$\mu|_{Y}$ to denote the restriction of this marking to $Y,$ by which 
we mean: 
\begin{itemize}
    \item If $\xi(Y)=-1$, then $\mu|_{Y}=\pi_{Y}(\mu)$.

    \item Otherwise, we let $\mu|_{Y}$ be the set of base curves which 
    meet $Y$ essentially, each taken with their associated transversal.
\end{itemize}

For a surface $Y$ with $\xi(Y)\geq 1$ and a simplex $v\subset\C(Y)$, 
we say that $X$ is a \emph{component domain}\index{component domain} 
of $(Y,v)$ if either: $X$ 
is a component of $Y\setminus v$ or $X$ is an annulus whose core curve is a 
component of $v$. More generally, we say that a subsurface $X\subset Y$ 
is a \emph{component domain of $g$} if for some simplex $v_{i}$ in $g$, $X$ 
is a component domain of $(D(g),v_{i})$. Notice that $v_{i}$ is 
uniquely determined by $g$ and $X$.

Furthermore, when $X$ is a component domain of $g$ and $\xi(X)\neq 0$, 
 define  the \emph{initial marking of $X$ relative to $g$}:
$$\I(X,g) = \left\{
\begin{array}{ll}
v_{i-1}|_{X} & v_i \mbox{ is not the first vertex (of $g$)} \\
\I(g)|_{X}& v_i \mbox{ is the first vertex} 
\end{array}
\right.
$$\index{$\I(X,g)$} 
Similarly define the \emph{terminal marking of $X$ relative to $g$} to be:
$$\T(X,g) = \left\{
\begin{array}{ll}
v_{i+1}|_{X} & v_i \mbox{ is not the last vertex} \\
\T(g)|_{X} & v_i \mbox{ is the last vertex} 
\end{array}
\right.
$$\index{$\T(X,g)$} 
Observe that these 
are each markings, since $\partial X$ is distance 1 in $\C(S)$ from 
each of $v_{i\pm 1}$, or in the case where $v_{i}$ is the first 
 vertex, then $\partial X$ is disjoint from 
$\base(\I(g))$ (similarly for the terminal 
markings)\index{marking!initial}\index{marking!terminal}.

\index{subordinate}
When $X$ is a component domain of $g$ with $\T(X,g)\neq\emptyset$, then 
we say that $X$ is \emph{directly forward subordinate} to $g$, 
written $X\fsubd g$\index{$X\fsubd g$}. Similarly, when $\I(X,g)\neq\emptyset$ we say that 
$X$ is \emph{directly backward subordinate} to $g$, 
written $g\bsubd X$\index{$g\bsubd X$}.

The definition generalizes to geodesics as follows.
\begin{defn} Let $g$ and $h$ be tight geodesics. We say that $g$ is {\em 
directly forward subordinate} to $h$, written $g\fsubd h$, when 
$D(g)\fsubd h$ and $\T(g)=\T(D(g), h)$. Similarly, $h$ is {\em 
directly backward subordinate} to $g$, written $g\bsubd h$, when 
$g\bsubd D(h)$ and $\I(h)=\I(D(h), g)$.
\end{defn}

We write \emph{forward subordinate}, or $\fsub$, to denote the transitive closure of 
$\fsubd$; similarly we define $\bsub$.

We can now define the main tool which was introduced in \cite{MasurMinsky:complex2}.

\begin{defn}
    A \emph{hierarchy (of geodesics)}  $H$, on $S$ is a collection of  tight 
    geodesics subject to the following constraints:\index{hierarchy}
    \begin{enumerate}
        \item There exists a tight geodesic whose support is $S$. This 
        geodesic is called the \emph{main geodesic}\index{main geodesic} 
	and is often denoted $g_{H}$\index{$g_{H}$}. 
	The initial and terminal markings of $g_{H}$ are 
        denoted $\I(H)$ and $\T(H)$.\index{marking!initial}\index{marking!terminal}
    
        \item Whenever there exists a pair of tight geodesics $g,k\in H$ and 
        a subsurface $Y\subset S$ such 
        that $g\bsubd Y \fsubd k$ then $H$ contains a unique tight geodesic $h$
	with domain $Y$ such that $g\bsubd h \fsubd k$.
    
        \item For each geodesic $h\in H$ other than the main 
        geodesic, there exists $g,k\in H$ so that  $g\bsubd h \fsubd k$.
    \end{enumerate}
\end{defn}

Using an inductive argument, it is proved in \cite{MasurMinsky:complex2} 
that given any two 
markings on a surface, there is a hierarchy which has initial marking 
one of them and terminal marking the other. 
The process begins with picking a base curve of the initial marking, one 
in the terminal marking, and a geodesic in $\C(S)$ between them (the 
main geodesic). Then, the second condition in the definition of a 
hierarchy forces certain proper subdomains to support geodesics: if 
there is a configuration $g \bsubd D \fsubd k$ such that $D$ does not support
any geodesic $h$ with $\I(h)$ = $\I(D,g)$ and $\T(h) = \T(D,k)$, then
we construct such a geodesic. When this geodesics is constructed, we 
can choose the initial vertex to be any element of $\base 
\I(Y,g)$ (similarly for the terminal vertex). This process continues until 
enough geodesic are included in $H$ so that the second and third 
conditions are satisfied. 
As the above suggests, there is in general 
not just one hierarchy, but many of them connecting any pair of 
markings. In our proof of 
\fullref{projest} (Projection estimates) we will exploit this flexibility by 
building our hierarchy subject to certain additional constraints which we find 
useful. 

When we consider several hierarchies at once, we use the 
notation $g_{H,Y}$\index{$g_{H,Y}$} to denote the geodesic of $H$ 
supported on $Y,$ when 
this geodesic exists it is  unique by \fullref{MMtwSST}.

For any domain $Y\subset S$ and hierarchy $H$, 
the \emph{backward} and \emph{forward sequences}
 are given respectively by $$\Sigma^{-}_{H}(Y)=\{b\in H: Y \subseteq D(b) \;\mbox{and}\;
 I(b)|_{Y}\neq \emptyset \}$$ $$\Sigma^{+}_{H}(Y)=\{f\in H: Y 
 \subseteq D(b) \; \mbox {and}\; T(f)|_{Y}\neq \emptyset \}.\leqno{\hbox{and}} $$
 The following theorem summarizes some results which are useful for 
 making computations with hierarchies. \index{$\Sigma^{\pm}_{H}(Y)$}

 \begin{thm}[Structure of Sigma \cite{MasurMinsky:complex2}]\label{MMtwSST}
 Let $H$ be a hierarchy, and $Y$ any domain in its support.

 \begin{enumerate}
 \item If $\Sigma^+_H(Y)$ is nonempty then it has the form of a
 sequence 
 $$f_0\fsubd\cdots\fsubd f_n=g_H,$$ 
 where $n\ge 0$. Similarly, 
 if $\Sigma^-_H(Y)$ is nonempty then it has the form of a
 sequence 
 $$g_H=b_m\bsubd\cdots\bsubd b_0,$$
 where $m\ge 0$.

 \item If 
 $\Sigma_{H}^\pm(Y)$ are both nonempty and $\xi(Y)\neq 0$, then $b_0 = f_0$, and
 $Y$ intersects every vertex of $f_0$ nontrivially.
 \item If $Y$ is a component domain of a geodesic $k\in H$ and
   $\xi(Y)\ne 0$, then
 $$f\in \Sigma_{H}^+(Y) \ \ \iff \ \ Y\fsub f,$$
 and similarly, 
 $$b\in \Sigma_{H}^-(Y) \ \ \iff \ \ b\bsub Y.$$
 If, furthermore, $\Sigma^\pm(Y)$ are both nonempty,
 then in fact $Y$ is the support of the geodesic $b_0=f_0$. 

 \item Geodesics in $H$ are determined by their supports. That is, if
   $D(h)=D(h')$ for $h,h'\in H$ then $h=h'$.

 \end{enumerate}
 \end{thm}

 Given a hierarchy, the following provides a useful criterion for determining 
 when a domain is the support of some geodesic in that hierarchy.

\begin {lem}[Large links \cite{MasurMinsky:complex2}]\label{MMtwLLL}
If $Y$ is any domain in $S$ and
$$d_{Y}(I(H),T(H)) >M_{2},$$ then $Y$ is the support of a geodesic $h$ in $H$, 
where $M_{2}$ only depends on 
the topological type of $S$.
\end{lem}

The geodesics in a hierarchy admit a partial ordering, which 
generalizes both the linear ordering on vertices in a geodesic and the 
ordering coming from forward and backwards subordinacy. Below we 
recall the basic definitions and a few properties of this ordering.

Given a geodesic $g$ in $\C_{1}(S)$ and a subsurface $Y,$ define the {\em footprint} 
of $Y$ in $g$, denoted $\phi_{g}(Y)$ to be the collection of vertices 
of $g$ which are disjoint from $Y.$  
It is easy to see that the diameter 
of this set is at most 2 and with a little work it can be shown that 
the footprint is always an interval (of diameter at most 2). The proof 
that $\phi_{g}(Y)$ is interval uses the assumption that the geodesic is 
tight: this result is the only place where the tightness assumption on 
geodesics is utilized. 
That footprints form intervals is useful as it allows one to make the 
following definition.

\begin{defn}\index{$g\prec_{t} k$}  For a pair of geodesics $g,k\in H$ we say {\em g 
precedes k in the time order}, or $g\prec_{t} k$, if there exists a 
geodesic $m\in H$ so that $D(g)$ and  $D(k)$ are both subsets of $D(m)$ and 
$$\max\phi_{m}(D(g))< \min\phi_{m}(D(k)).$$
\noindent
We call $m$ the \emph{comparison geodesic}\index{comparison geodesic}.
\end{defn}

Time ordering is a (strict) partial ordering on geodesics in a hierarchy, 
this is proven in \cite{MasurMinsky:complex2}. It is worth remarking that 
when a pair of geodesics are time ordered, they are time ordered 
with respect to a  unique comparison geodesic.

The following provides a way to use the time ordering to gain 
information about the hierarchy.

\begin {lem}[Order and projections 
\cite{MasurMinsky:complex2}]\label{MMtwOPL}
Let $H$ be a hierarchy in $S$ 
and $h,k \in H$ with $D(h)=Y$ and $D(k)=Z$.  Suppose that $Y \cap Z
\neq \emptyset$ and neither domain is contained in the other.  Under 
these conditions, if
$k \prec_{t} h$ then $d_{Y}(\partial Z, I(H)) \leq M_{1}+2$ and
$d_{Z}(T(H), \partial Y) \leq M_{1}+2$. The constant $M_{1}$ only depends on 
the topological type of $S$.
\end{lem}

Since we will often use the above hypothesis, we introduce the terminology 
that a pair of subsurfaces $Y$ and $Z$ of $S$ \emph{overlap}\index{overlap} when 
$Y \cap Z \neq \emptyset$ and neither domain is contained in the other.
A application of this lemma which we shall often use is: given 
$h_{1},h_{2},k \in H$ with $D(h_{i})=Z_{i}$ and $D(k)=Y$ and for which 
$h_{i}\prec_{t} k$, then $d_{Y}(\partial Z_{1},\partial Z_{2})\leq 
2(M_{1}+2)$.

The next result provides a way to translate distance computations in the 
mapping class group to computations in curve complexes of subsurfaces.

\begin{thm}[Move distance and projections 
    \cite{MasurMinsky:complex2}]
    \label{mcgalmostquasi}
There exists a constant $t(S)$ such that for each $\mu, \nu \in \M$ and any
threshold $t>t(S)$, there exists $K(=K(t))$ and $C(=C(t))$ such that:
$$\frac{1}{K}d_{\M}(\mu,\nu) -C
\leq \sum_{\substack{Y\subseteq S\\ d_{Y}(\mu,\nu)>t}} d_{Y}(\mu,\nu)
\leq K d_{\M}(\mu,\nu) +C.$$
\end{thm}

The importance of \fullref{mcgalmostquasi} can not be understated, 
as this theorem provides the crucial result that hierarchies give rise 
to quasi-geodesic paths in the mapping class group. We now indicate  
this construction.

Given a hierarchy $H$, \cite{MasurMinsky:complex2} define a 
\emph{slice}\footnote{This is called 
a \emph{complete slice} in \cite{MasurMinsky:complex2}.} 
in $H$ to be a set $\tau$ of pairs $(h,v)$, where $h\in H$ and $v$ is a 
vertex of $h$, satisfying the following conditions:
\begin{enumerate}
    \item A geodesic $h$ appears in at most one pair of  $\tau$.

    \item There is a distinguished pair $(h_{H},v_{\tau})$ 
	appearing in $\tau$, where $h_{H}$ is the main geodesics of $H$.

    \item For every $(k,w)\in\tau$ other than $(h_{H},v_{\tau})$, 
     $D(k)$ is a component domain of $(D(h),v)$ for some $(h,v)\in\tau$.

    \item Given $(h,v)\in\tau$, for every component domain $Y$ of 
    $(D(h),v)$ there is a pair  $(k,w)\in\tau$ with $D(k)=Y$.
\end{enumerate}

To any slice $\tau$ is associated a unique marking as follows. 
All the vertices which appear in geodesics of non-annular subsurfaces 
are distinct and disjoint, so these form a maximal simplex of $\C(S)$ 
which we take to be $\base(\mu_{\tau})$. Then, by the forth condition, 
$\tau$ must contain a pair of the form $(h_{\gamma},t)$ 
for each $\gamma\in\base(\mu_{\tau})$. Letting $t$ be the transversal 
for $\gamma$ we obtain a complete marking, which we denote  $\mu_{\tau}$. 
This marking may not be clean, but a lemma of \cite{MasurMinsky:complex2} 
implies that a complete clean marking can be found in a uniformly bounded 
neighborhood; thus we assume  $\mu_{\tau}\in\M(S)$.

The set of slices admits a partial ordering (relating to the time 
ordering discussed above), which then descends to provide a partial 
ordering on the markings $\mu_{\tau}$ for a given hierarchy $H$. 
In  \cite{MasurMinsky:complex2}  it is  
shown that slices nearby in the slice ordering yield markings which are 
close in $\M(S)$; combined with \ref{mcgalmostquasi} these sequences 
of markings form quasi-geodesic paths in $\M(S)$, called 
\emph{hierarchy paths}.

\subsection{Teichm\"{u}ller space and the pants complex}

For a detailed reference on Teichm\"{u}ller space, consult
Imayoshi and Taniguchi \cite{ImayoshiTaniguchi}.

\begin{defn} For a fixed open topological surface $S$, the
    \emph{Teichm\"{u}ller space of $S$}\index{Teichm\"{u}ller space} 
    is the space of equivalence 
    classes of pairs $(X,f)$, where $X$ is a finite area hyperbolic surface and
    $f\co S\to X$ is a homeomorphism. 
    A pair $(X_{1},f_{1})$ and 
    $(X_{2},f_{2})$ are considered equivalent if there exists an 
    isometry $h\co X_{1}\to X_{2}$ such that 
    $h\circ f_{2}$ is homotopic to $f_{1}$.
\end{defn}

A topology is obtained on Teichm\"{u}ller space by infimizing over the 
distortion of maps $h\co X_{1}\to X_{2}$.
Topologically Teichm\"{u}ller space is fairly easily understood, as 
the following classical result indicates:

\begin{thm} Teichm\"{u}ller space is homeomorphic to $\R^{6g-6+2p}$.
\end{thm}

There are several natural metrics on Teichm\"{u}ller space and its
metric structure is far less transparent than its topological
structure. Often of interest is the natural complex structure which
Teichm\"{u}ller space carries.  What will be important in the sequel
is the \emph{Weil--Petersson metric}\index{Weil--Petersson metric} on
Teichm\"{u}ller space, which is a K\"{a}hler metric with negative
sectional curvature. We will not need the integral form of the
definition, so we just mention that the metric is obtained by
considering a natural identification between the complex cotangent
space of Teichm\"{u}ller space and the space of holomorphic quadratic
differentials, then defining the Weil--Petersson metric to be the one
dual to the $L^{2}$--inner product on the space of holomorphic
quadratic differentials. See Wolpert \cite{Wolpert:survey} for a
survey on the Weil--Petersson metric and its completion.

\begin{defn}
A \emph{pair of pants} is a thrice punctured sphere. 
A \emph{pants decomposition}\index{pants decomposition}  
of a surface $S_{g,p}$ is a maximal 
collection of pairwise non-parallel 
homotopy classes of simple closed curves; this decomposition obtains its name from 
the observation that such a curve system 
cuts  $S_{g,p}$ into $2g+p-2$ pairs of pants.
\end{defn}

It is easy to verify that any pants decomposition consists of exactly 
 $3g+p-3$ disjoint simple closed curves on $S$ and is thus a maximal 
 simplex in $\C(S)$. These decompositions have long been 
useful in the study of the mapping class group and Teichm\"{u}ller space.

\begin{figure}[ht!]\label{figurepairofpants}
\centerline{\includegraphics[width=2.5truein]{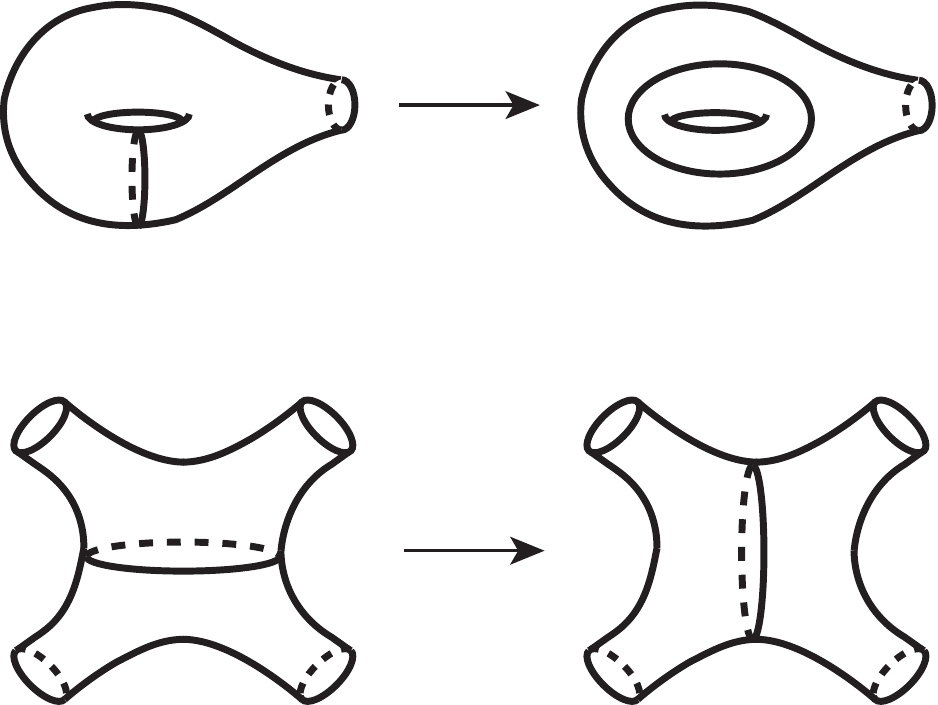}}
\caption{Elementary moves in the pants complex}
\end{figure}

Originally defined by Hatcher \cite{Hatcher:Pants}, the 
 \emph{pants complex of $S$}\index{pants complex}, denoted 
 $\P(S)$\index{$\P(S)$}, 
is a way of comparing all possible pants decompositions on a fixed 
surface. This complex consists of a vertex for every pants decomposition 
and an 
edge  between each pair of decompositions which differ by an 
\emph{elementary move}\index{pants complex!elementary move}. Two pairs of pants 
$P=\{\gamma_{1},\ldots,\gamma_{3g+p-3}\}$ and 
$P'=\{\gamma'_{1},\ldots,\gamma'_{3g+p-3}\}$ differ by an
\emph{elementary move}
if $P$ and $P'$ can be reindexed so that both:
\begin{enumerate}
    \item $\gamma_{i}=\gamma'_{i}$ for all $2\leq i\leq 3g+p-3$

    \item In the component $T$ of $S\setminus \cup_{2\leq i\leq 
    3g+p-3}\gamma_{i}$ which is not a thrice punctured sphere 
    ($T$ is necessarily either a once punctured torus or a 
    4 punctured sphere) we have (see \fullref{figurepairofpants}):
    $$d_{\C(T)}(\gamma_{1},\gamma'_{1})=1.$$
\end{enumerate}
\noindent

$\P(S)$ is metrized  by giving each edge the metric of the Euclidean 
interval $[0,1]$.
Our interest in this space comes out of the following remarkable 
theorem of Brock:
\begin{thm}{\rm\cite{Brock:wp}}\qua $\P(S)$ is quasi-isometric to the Teichm\"{u}ller space 
of $S$ with the Weil--Petersson metric.
\end{thm}

\begin{rmk}\label{pantsvsmcg} 
    Noticing that  a marking without transverse data is 
    just a pants decomposition,
    in Section 8 of \cite{MasurMinsky:complex2} it is remarked that 
    all of the constructions in their paper using hierarchies 
    to obtain results about the marking complex can be replaced 
    with analogous theorems for the pants complex. 
    This is done by replacing markings (and hierarchies, etc)
    by markings without transverse data (and hierarchies-without-annuli, 
    etc). We will often implicitly use this result to note that
    results we prove from the mapping class group also hold in Teichm\"{u}ller 
    space (in particular, this allows one to obtain Teichm\"{u}ller 
    space analogues of the results in 
    Sections~\ref{sectionhyperbolicitylowcomplexity} and~\ref{sectiontransversaltree}). 
    The following is the main result of this form which we will
    explicitly use.
\end{rmk}

\begin{thm}[Move distance and projections for $\P(S)$
    \cite{MasurMinsky:complex2}] \label{pantsalmostquasi}
There exists a constant $t(S)$ such that for each $\mu, \nu \in \P$ (and any
threshold $t>t(S)$) there exists $K(=K(t))$ and $C(=C(t))$ such that:
$$\frac{1}{K}d_{\P}(\mu,\nu) -C
\sum_{ \substack{\text{non-annular } Y \subseteq S \\ d_{Y}(\mu,\nu)>t}}
d_{Y}(\mu,\nu)
\leq K d_{\P}(\mu,\nu) +C.$$
\end{thm}

\section{Projection estimates}\label{chapterProjectionEstimates}

In this section we address the following question: given a pair 
of subsurfaces $X_{1},X_{2}$ of $S$, what is the 
image of $\M(S)$ under the map which projects markings into the 
product of curve 
complexes of these subsurfaces, ie, the image of $\M(S)$ as a subset of 
$\C(X_{1})\times\C(X_{2})$?  
We answer this in Theorems~\ref{projest} and 
\ref{geomprojest} using hierarchies in the curve complex as 
our main tool (see \fullref{ccbackground} for background). 
The geometric description 
underlying \fullref{geomprojest} can be generalized 
to arbitrary finite 
collections of subsurfaces by introducing an notion called a 
coarse tree-flat and using a delicate inductive step, this analysis 
is carried out in \cite{Behrstock:thesis}. 

We begin this section with a pair of toy examples which 
illustrate the dichotomy of 
 \fullref{projest} (Projection estimates), but which can be proven without 
use of the full hierarchy machinery. In \fullref{sectionprojest} 
 we prove \fullref{projest} and its 
geometric analogue \fullref{geomprojest}.

\subsection{A motivating example: the dichotomy}
\label{motivatingexamples}

Here we give two computations on the genus three surface. We exhibit a 
dichotomy between the behavior of projection maps into  
disjoint and intersecting pairs of subsurfaces. These examples  
illustrate the geometric meaning of \fullref{projest} 
and provide a useful warm up for its proof, 
as the arguments we use here  are a particularly easy form 
of what in general is a technical argument involving hierarchies.

\begin{exa}[Disjoint subsurfaces] \label{exampledisjoint}
    Let $S$ be a closed genus three surface and let $X$ and $Y$ be 
    subsurfaces of $S$ which are each once punctured tori 
    and which can be realized disjointly on $S$ (\fullref{disjointsurfaces}).

    \begin{figure}[ht!]\label{disjointsurfaces}
    \labellist \hair 5pt
    \pinlabel {$X$} at 49 130
    \pinlabel {$Y$} at 235 29
    \pinlabel {$\partial X$} at 95 59
    \pinlabel {$\partial Y$} at 191 100
    \pinlabel {$\overbrace{\hspace*{1in}}$} [t] at 49 130
    \pinlabel {$\underbrace{\hspace*{1in}}$} [b] at 235 29
    \endlabellist
    \centerline{\includegraphics[width=3truein]{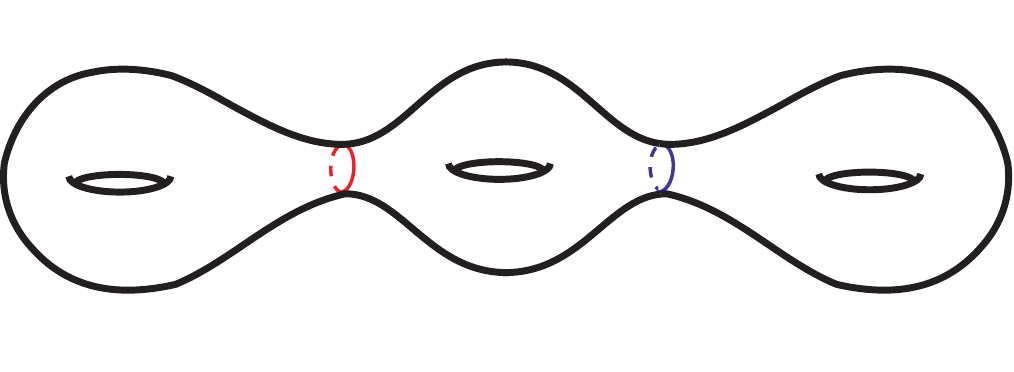}}
    \caption{$S$ and two disjoint subsurfaces}
    \end{figure}
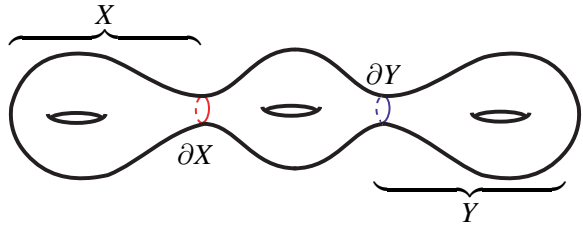
    
    Let $\alpha\in\C(X)$ and $\beta\in\C(Y)$. Note that both these 
    curves can also be considered as elements of $\C(S)$, then 
    since $X\cap Y=\emptyset$ 
    it follows that $d_{\C(S)}(\alpha,\beta)=1$. Thus $\alpha$ and 
    $\beta$ are vertices of a common top dimensional simplex $\mu\subset\C(S)$. 
    Taking a complete  clean marking $\mu'\in\M(S)$ with 
    $\base(\mu')=\mu$, we have that 
    the map $$\M(S) \to \C(X)\times \C(Y)$$ is onto.

\end{exa}

The next example provides the more interesting half of the dichotomy.

\begin{exa}[Overlapping subsurfaces]\label{exampleoverlap}
    Let $S$ be the closed genus three surface and consider two 
    subsurfaces $X$ and $Y$ which each have genus two and one puncture.
   Moreover, suppose that $X$ and $Y$ overlap in a twice punctured 
   torus, $Z$.
 Letting  $\mu$ be an element of $\C(S)$ which intersects both $X$ and 
 $Y$. We shall show that both $d_{X}(\mu,\partial Y)$ and
$d_{Y}(\mu,\partial X)$ cannot simultaneously be large.

  \begin{figure}[ht!]\label{overlappingsurfaces}
    \labellist \hair 5pt
    \pinlabel {$X$} at 97 130
    \pinlabel {$Y$} at 189 27
    \pinlabel {$\partial X$} [b] at 192 88
    \pinlabel {$\partial Y$} [t] at 99 72
    \pinlabel {$\overbrace{\hspace*{2in}}$} [t] at 97 130
    \pinlabel {$\underbrace{\hspace*{2in}}$} [b] at 189 27
    \endlabellist
  \centerline{\includegraphics[width=3truein]{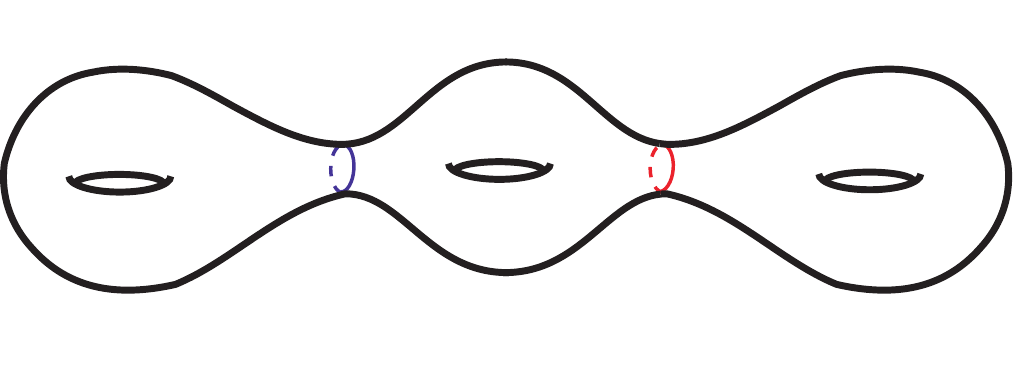}}
  \caption{$S$ and two overlapping subsurfaces}
\end{figure}
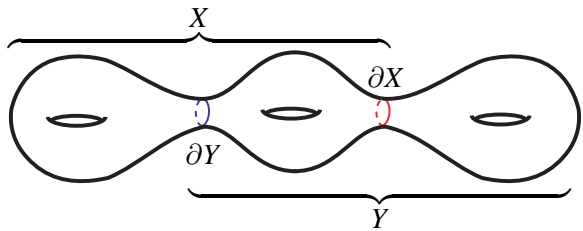

Let $g$ be a geodesic in $\C_{1}(S)$ connecting $\mu$ to $\partial X$, with vertices
$v_{0}, v_{1}, \ldots, v_{n}$ where $v_{0}=\mu$ and $v_{n}=\partial X$. 
Since $g$ is a geodesic in $\C(S)$, the curve $v_{n-1}$ may be 
disjoint from $X$, 
but for all $k<n-1$ we have
$v_{k}\cap \partial X\neq \emptyset$ and in particular, 
$\pi_{X}(v_{k})\neq\emptyset$ and $\pi_{Y}(v_{k})\neq\emptyset$. 
Since $S=X\cup Y$ we have either 
$v_{n-1}\cap Y \neq\emptyset$ or $v_{n-1}\subset X$.

If $v_{n-1} \cap Y \neq\emptyset$, then every vertex of $g$ intersects $Y$ 
non-trivially. \fullref{MMtwbgi} (Bounded geodesic image 
\cite{MasurMinsky:complex2}) states that if $g$ is a geodesic in $\C_{1}(S)$ 
with every vertex non-trivially intersecting a subsurface $Y\subset S$, then 
$\diam_{Y}(\pi_{Y}(g))<K$ for a constant $K$ depending only on the 
topological type of $S$. Accordingly, this theorem implies that
we have 
$d_{Y}(\mu,\partial X)<K$ for a constant $K$ depending
only on the genus of $S$.

If $v_{n-1}\cap Y =\emptyset$, then it must be the case 
that $v_{n-1}\subset X$. In this case, we have $v_{k}\cap 
X\neq\emptyset$ for each $0\leq k\leq n-1$.  Applying the Bounded
Geodesic Image Theorem tells us that
$d_{X}(\mu,v_{n-1})<K$.  Since $v_{n-1}\cap Y=\emptyset$, we have in 
particular that  $v_{n-1}\cap \partial Y=\emptyset$ which implies that 
$d_{X}(v_{n-1},\partial Y)\leq 3$.  Thus 
we may conclude that $d_{X}(\mu,\partial Y)<K+3$.

Putting the two cases together, this example shows that there exists 
a constant $C=K+3$,
for which any $\mu\in\C(S)$ must satisfy either 
$ d_{X}(\mu,\partial Y)<C$ or $d_{Y}(\mu,\partial X)<C.$

If we take $\mu$ to be one of the 
base curves for some complete clean marking $\nu$ of $S$, we have shown 
that for any marking $\nu\in\M(S)$ either 
$$ d_{X}(\nu,\partial Y)<C \mbox{ or } d_{Y}(\nu,\partial X)<C.$$

\end{exa}

We consider this a toy example, since the only tool we needed was
the Bounded geodesic image Theorem.  
Difficulties in generalizing this example to arbitrary 
overlapping subsurfaces include problems such as:
 the boundaries of these subsurfaces may fill $S$ or the 
union of these subsurfaces may not be all of $S$---both properties that 
were crucial to the above argument. To deal with these problems we  
use a powerful technical tool called a hierarchy, in which  
one  considers not just a geodesic in the curve complex of $S$, but a family 
of geodesics each in the curve complex of some subsurface of $S$. 
The tools introduced by Masur and Minsky for calculating with these 
allow one to give the argument in the next section
which is obtained by bootstrapping the key idea used in the above example.

\subsection{Projection estimates theorem}\label{sectionprojest}

\begin{thm}[Projection estimates]\label{projest}
    Let $Y$ and $Z$ be two overlapping surfaces in $S$,
    with  $\xi(Y)\neq 0 \neq \xi(Z)$, 
    then for any $\mu\in\M(S)$:
    $$d_{Y}(\partial Z, \mu)>M  \implies d_{Z}(\partial Y, 
    \mu)\leq M,$$
    for a constant $M$ depending only on the topological type of $S$.
\end{thm}

\begin{proof} Let $M_{1}$ and $M_{2}$ be given by
    Lemmas~\ref{MMtwLLL} and \ref{MMtwOPL} 
     (Large link Lemma and Order and projections Lemma 
     \cite{MasurMinsky:complex2}). Define: $$M=\max\{M_{1}+2,M_{2}+3\}.$$
    We now fix a point $\mu\in\M$ for which $d_{Y}(\mu, \partial Z) > M$ 
    and will show that this 
    implies a bound on  $d_{Z}(\mu, \partial Y)$.

    Consider a hierarchy $H$ with  
    $\base(I(H))\supset\pi_{Z}(\mu) \cup \partial Z$ and $T(H)=\mu$
    (when $Z$ is an  annulus take 
    $\base(I(H))=\partial Z$ with transversal $\pi_{Z}(\mu)$).
    Furthermore, we assume  $H$ is built subject to the following two constraints. 
    (See \fullref{hierarchiesbackground} for a discussion of the 
    choices involved when building a hierarchy.)

    \begin{quote}
        \begin{itemize}
            \item Choose the initial vertex of $g_{H}$ to be an element 
            of $\partial Z$ which intersects $Y$. Call this vertex $v_{0}$.
        
            \item When adding a geodesic to $H$ with a given initial 
            marking $I(Y,g)$, if 
            $v_{0}\in\base(I(Y,g))$ choose this as the initial vertex 
            for the geodesic. If $v_{0}\notin\base(I(Y,g))$, 
	    then if any elements of $\partial Z$ are in 
             $\base(I(Y,g))$ choose one of these to be the initial vertex 
             for the geodesic.
	\end{itemize}
    \end{quote}

 The proof is now broken into four steps. We show 
 there exists a geodesic $k\in H$ with domain $Z$ and a geodesic $h\in H$ with 
 domain $Y$. Then we show $k \prec_{t} h$, from which the theorem 
 follows as a consequence of  \fullref{MMtwOPL} 
 (Order and projections \cite{MasurMinsky:complex2}).

\medskip
{\bf Step (i)}\qua {\sl There exists a geodesic $k$ with domain $Z$.}
  \smallskip

 We start by considering the forward and backward  sequences, 
$\Sigma^{+}_{H}(Z)$ and $\Sigma^{-}_{H}(Z)$, as defined in the  
discussion preceding \fullref{MMtwSST} 
(Structure of Sigma \cite{MasurMinsky:complex2}).
 
 Since $\mu$ is complete we have $\pi_{Z}(\mu) \neq \emptyset$; thus both 
 the initial and terminal markings restrict to give nontrivial markings 
 on $Z$ (when $Z$ is an annulus, we are using that $I(H)$ has a 
 transverse curve). 
 This implies for any subsurface $Q\supseteq Z$ that
 $\Sigma^{-}_{H}(Q)$ and $\Sigma^{+}_{H}(Q)$
 both contain $g_{H}$ and thus in particular are
 nonempty. \fullref{MMtwSST}  shows that when
 $Q$ is a component domain of a geodesic in $H$ and both
 $\Sigma^{+}_{H}(Q)$ and $\Sigma^{-}_{H}(Q)$ are non-empty, then $Q$
 must be the support of a geodesic in $H$.

 Since the first vertex of $g_{H}$ is $v_{0}\in\partial Z$, we know it 
 has  
 a component domain $Q_{1} \subsetneq S$ which contains $Z$. 
 By the above observation, we know 
 that $\Sigma^{+}_{H}(Q_{1})$ and $\Sigma^{-}_{H}(Q_{1})$ are both 
 non-empty, and thus $Q_{1}$ supports a geodesic which we call $k_{1}$. 
 If $Q_{1}=Z$ we have produced a geodesic supported on $Z$ and are done, 
 otherwise since $Q_{1}\supsetneq Z$ we can choose  an 
 element of $\partial Z$ as the 
 first vertex of the geodesic supported in $Q_{1}$.

 Starting from the base case $Q_{0}=S$ and $k_{0}=g_{H}$, the above 
  argument produces 
 a sequence of properly nesting subsurfaces $S=Q_{0}\supsetneq 
 Q_{1}\supsetneq \ldots\supsetneq Q_{n}=Z$ where for each $i>0$ the 
 subsurface $Q_{i}$ is
 a component domain of a geodesic $k_{i-1}\in H$
 supported in $Q_{i-1}$ and each of the $k_{i}$ has an element of $\partial 
 Z$ as an initial vertex.

 Since $S$ is a surface of finite type, the above sequence of nested 
 surfaces must terminate with $Z$ after finitely many steps.
 
 Thus $Z$ is a component domain of a geodesic in $S$ and supports a 
 geodesic which we call $k$.

  \medskip
{\bf Step (ii)}\qua {\sl There exists a geodesic $h$ with domain $Y$.}
  \smallskip
 
  By hypothesis we have $d_{Y}(\partial Z, \mu) >M_{2}+3$ and thus:
 $$d_{Y}(I(H), T(H))=d_{Y}( \pi_{Z}(\mu) \cup \partial Z, \mu) 
 >M_{2}.$$ 
 Together with 
  \fullref{MMtwLLL} (Large link Lemma \cite{MasurMinsky:complex2}), this 
 implies $Y$ is the support of a geodesic $h$ in $H$.

 \medskip
 \noindent
 {\bf Step (iii)}\qua {\sl $k \prec_{t} h$}
  \smallskip
 
 In this step we will prove that the geodesic $k$ precedes $h$ in the 
 time ordering on geodesics in $H$.

 Since $D(k)=Z$ and $D(h)=Y$ and each are contained in $D(g_{H})=S$:
 if 
 $\max \phi_{g_{H}}(Z)$ $< \min \phi_{g_{H}}(Y)$, then 
 $k \prec_{t} h$ which is what we wanted to prove. In the general case, we 
 provide an inductive procedure to  show that these geodesics have 
 the desired time ordering.
 
 We refer to the ordered vertices of $g$ as $v_{i}(g)$. Recall  
 that by the first part of the constraint we have 
 $v_{0}(g_{H})=v_{0}$,  
 where $v_{0}$ was chosen to satisfy 
  $v_{0} \in \partial Z$ and $v_{0}\cap 
 Y\neq\emptyset$.  Also, since
 $g_{H}$ is a tight geodesic, $v_{1}(g_{H})=\partial
 F(v_{0}(g_{H}),v_{2}(g_{H}))$ where $F(\alpha,\beta)$ denotes the surface filled by
 $\alpha$ and $\beta$.
 
 Summarizing, we have:
 
\begin{itemize}
 \item $Y \cap v_{0} \neq \emptyset$, and thus $v_{0} \notin \phi_{g_{H}}(Y)$.
 \item $Z$ is contained in a component domain of $v_{0}$, thus $v_{0} \in
 \phi_{g_{H}}(Z)$.
 \item Since $g_{H}$ is a geodesic, we know $v_{2}(g_{H})$
 must intersect $v_{0}(\in \partial Z)$. Since the diameter of a 
 footprint is at most 2, it now follows that 
  $\max \phi_{g_{H}}(Z)\leq v_{1}(g_{H})$.
\end{itemize}
 
Together these imply: $$\max\phi_{g_{H}}(Z)\leq v_{1}(g_{H})\mbox{ and }
 v_{1}(g_{H})\leq\min \phi_{g_{H}}(Y).$$
Now there are two mutually exclusive cases to consider:
\begin{enumerate}
    \item \label{middleofproofcaseon}
    $v_{1}(g_{H}) \cap Z \neq \emptyset$ or $v_{1}(g_{H})\cap Y\neq \emptyset$

    \item \label{middleofproofcasetw}
    $v_{1}(g_{H}) \cap Z=\emptyset=v_{1}(g_{H})\cap Y$.
\end{enumerate}

In the first case, depending on which of the two sets is non-empty, we 
obtain  $v_{1}(g_{H}) 
\notin\phi_{g_{H}}(Z)$ or $v_{1}(g_{H})\notin\phi_{g_{H}}(Y)$, 
respectively. Either of these imply  
$\max \phi_{g_{H}}(Z) < \min \phi_{g_{H}}(Y)$, proving that in 
case~\ref{middleofproofcaseon} we have $k \prec_{t} h$.

	Observing that in $S_{1,1}$ and $S_{0,4}$ the footprint of a domain must 
	consist of at most one vertex, we have  
	$v_{1}(g_{H}) \notin\phi_{g_{H}}(Z)$ and thus when $S$ is 
	either of these surfaces we are in case~\ref{middleofproofcaseon} 
	and thus $k \prec_{t} h$. 
	As the annulus does not contain any pair of surfaces which overlap we
	may now assume for the rest of the proof that $\xi(S)>1$. In 
	particular, for the remainder of the argument we need not consider the 
	sporadic cases where $\C(S)$ has a special 
	definition and thus we can assume two distinct homotopy 
	classes of curves have
	distance one in the curve complex if and only if they can be realized 
	disjointly.

	For the remainder of this step we  
	assume that $v_{1}(g_{H})\in\phi_{g_{H}}(Z)$ and 
	$v_{1}(g_{H})\in\phi_{g_{H}}(Y)$ and will prove that we still obtain $k \prec_{t} h$.
	Since $v_{1}(g_{H})\cap Z=\emptyset$ there is a
	component of $S \setminus v_{1}(g_{H})$ which contains $Z$ and since $Z
	\cap Y \neq \emptyset$ the same component contains $Y$ as well---we will call
	this component $W_{1}$.  Since $W_{1}$ is a component domain of $H$ 
	which intersects both $I(H)$ and $T(H)$ (since $W_{1}\supset Z$), 
	\fullref{MMtwSST} (Structure 
	of Sigma) implies that it supports a geodesic $l_{1}$.
	
	Since $W_{1}$ is a component domain of $v_{1}(g_{H})$, we have 
	$v_{0}\in I(W_{1},g_{H})$. Thus, by our convention for choosing geodesics 
	in $H$, we choose $v_{0}(l_{1})=v_{0}$.
	As before we have
	$v_{0}(l_{1})\in \phi_{l_{1}}(Z)$, $v_{2}(l_{1}) \notin
	\phi_{l_{1}}(Z)$, and $v_{0}(l_{1}) \notin \phi_{l_{1}}(Y)$.  If 
	$v_{1}(l_{1})$ is
	in both $\phi_{l_{1}}(Z)$ and $\phi_{l_{1}}(Y)$ then we again restrict
	ourselves to the appropriate component $W_{2}$ of $W_{1} \setminus
	v_{1}(l_{1})$ and repeat the argument with $W_{2}$ and $l_{2}$.  This
	gives a properly nested collection of subsurfaces each of which
	contains $Y \cup Z$, so the process terminates in a finite number of
	steps to produce a 
	geodesic $l_{n}$ with domain $W_n \supseteq Y\cup Z$ and which satisfies
	$v_{0}(l_{n})=v_{0} \in \phi_{l_{n}}(Z)$, 
	$v_{2}(l_{n})\notin \phi_{l_{n}}(Z)$, 
	$v_{0}(l_{n}) \notin \phi_{l_{n}}(Y)$, 
	and either $v_{1}(l_{n})\notin \phi_{l_{n}}(Z)$
	or $v_{1}(l_{n}) \notin \phi_{l_{n}}(Y)$. Hence $\max \phi_{l_{n}}(Z) < \min 
	\phi_{l_{n}}(Y)$, and we have $k\prec_{t} h$ with comparison geodesic $l_{n}$.

   \medskip
\textbf{Step (iv)}\qua {\sl Conclusion}
  \smallskip
  
 We have now produced a hierarchy $H$ with 
 geodesics $k$ and $h$ with $D(h)=Y$, $D(k)=Z$, 
 and $k\prec_{t} h$. Thus  
 \fullref{MMtwOPL} (Order and projections \cite{MasurMinsky:complex2})
 implies that $d_Z (\partial Y, T(H)) \leq M_{1}+2$.

 Since $T(H)=\mu$, this implies
 $d_Z(\partial Y, \mu) \leq M_{1}+2$ which is exactly what we wanted
 to show.
    \end{proof}

It would be interesting to have a more elementary proof of 
\fullref{projest}, although the given argument is useful in that 
it provides the reader 
an introduction to the tools we need for the later parts
of this work.  C~Leininger has suggested to me ideas for a more
constructive proof.

\subsubsection{Geometry of projection estimates}\label{sectiongeomprojest}

The main case of the following result is an immediate corollary of 
 \fullref{projest} (Projection estimates); we give this restatement 
because it emphasizes the geometry underlying that result.

    \begin{thm}[Projection estimates, geometric version]\label{geomprojest}
	For any distinct 
    subsurfaces $Y$ and $Z$ of 
$S$, exactly one of the following holds for the map $$\pi_{Y} \times
\pi_{Z}\co \M(S) \to \C(Y) \times \C(Z)$$
 \begin{enumerate}
	\item $Y \cap Z =\emptyset$ in which case the map is onto. 
	\item One of the surfaces is contained inside the other, say $Z 
	\supset Y$. Then  
	the image is contained in a radius 3 neighborhood of the set
	$$\C(Y) \times \pi_{Z}(\partial Y) \: \bigcup \:
	\Graph(\pi_{Z\to Y}),$$ where $\pi_{Z\to Y}$ is the
	projection map from $\C(Z)\setminus B_{1}(\partial Y)$ to $\C(Y)$.
	Furthermore, $\pi_{Y} \times\pi_{Z}(\M(S))$ is $6$--dense in this set.
	\item $Y$ and $Z$ overlap, in which case 
	the image is contained in a radius $M=\max\{M_{1}+2,M_{2}\}$
	neighborhood of the set
	$$\pi_{Y}(\partial Z) \times \C(Z) \:\bigcup\:
	\C(Y) \times \pi_{Z}(\partial Y).$$ The constant $M$ 
	depends only on the topological type of $S$.
	Furthermore, $\pi_{Y} \times\pi_{Z}(\M(S))$ is $6\cdot M$--dense in 
	this set.
    \end{enumerate}
\end{thm}
 
\begin{proof}\ \
    
   \noindent{\bf Case 1\qua $Y\cap Z=\emptyset$}  

    Any pair consisting of a curve in $\C(Z)$ and a curve in 
    $\C(Y)$ can be completed to a complete marking on $S$ so the map 
    is onto.

    \medskip
    \noindent{\bf Case 2\qua $Y\subset Z$}

    For any $\mu\in\M$ we have that $d_{Z}(\pi_{Z}'(\mu),\pi_{Z}(\mu)) \leq 1$ 
    and thus when we restrict these arcs (curves) to $Y$ we get that 
    $d_{Y}(\pi_{Y}'\pi_{Z}'(\mu),\pi_{Y}' \pi_{Z}(\mu)) \leq 1$.  
    Applying \fullref{MMtwlpl} (Lipschitz projection
    \cite{MasurMinsky:complex2}), 
    we have  
     $d_{Y}(\phi \pi_{Y}' \pi_{Z}'(\mu),\phi \pi_{Y}' \pi_{Z}(\mu))
    \leq 3$. Thus we have $d_{Y}(\mu,\pi_{Z\to Y}\circ\pi_{Z}(\mu))\leq 3.$

This  proves that if $Y\subset Z$ then 
$\pi_{Y} \times \pi_{Z}(\M)$ is contained in a radius 3 neighborhood of 
$\C(Y) \times \pi_{Z}(\partial Y) \: \bigcup \: \Graph(\pi_{Z\to Y})$.

We now prove the density statement. Let $(\gamma,\delta)\in \Graph(\pi_{Z\to 
Y})\subset \C(Y)\times\C(Z)$. Then take any $\mu\in\M(S)$ 
which has $\delta\in\base(\mu)$. By \fullref{MMtwlpl} 
we have $\pi_{Z}(\mu)\subset  B_{3}(\delta)$. By hypothesis we have 
$\gamma\in\pi_{Y}(\delta)$, so  \fullref{MMtwlpl} also shows 
$\pi_{Y}(\mu)\subset B_{3}(\gamma)$. Thus we have 
$(\gamma,\delta)\subset B_{6}(\pi_{Y}(\mu) \times \pi_{Z}(\mu))$. Similarly, 
given any $\gamma\in\C(Y)$, since $Y\subset Z$ we have $\gamma$ and 
$\partial Y$ are disjoint. Thus if 
$\mu\in\M(S)$ is chosen so that both $\gamma$ and 
$\partial Y$ are in $\base(\mu)$, we 
have $\pi_{Z}(\mu)\subset B_{3}(\gamma)$ and 
$\pi_{Y}(\mu)\subset B_{3}(\gamma)$. Thus $(\gamma,\pi_{Z}(\partial 
Y))\subset B_{6}(\pi_{Y}(\mu) \times \pi_{Z}(\mu))$. Thus we have shown that 
$\C(Y) \times \pi_{Z}(\partial Y) \: \bigcup \:	\Graph(\pi_{Z\to Y}) 
\subset  B_{6}(\pi_{Y} \times \pi_{Z}(\M))$.

    \medskip
    \noindent{\bf Case 3}\qua $Y\cap Z$, $Y\nsubseteq Z$, and $Z\nsubseteq Y$
   
    If $\mu\in\M(S)$ projects to $\C(Y)$ with
    $d_{Y}(\partial Z, \mu)>M$ then 
    \fullref{projest} (Projection estimates)  
    implies that $d_{Z}(\partial Y, \mu)\leq M$. Thus the image of the map 
    $\M(S)\to\C(Y)\times \C(Z)$ is contained in the union of the 
    radius $M$ neighborhood of $\pi_{Y}(\partial Z)\times \C(Z)$ 
    with the radius $M$ neighborhood of $\C(Y)\times\pi_{Z}(\partial Y)$ 
    as claimed.

    Given any $\gamma\in\C(Y)$, \fullref{MMtwlpl} implies that 
    $\diam_{Z}(\pi_{Z}(\gamma), \pi_{Z}(\partial Y))\leq 3$. If 
    $\mu\in\M(S)$ is chosen so that both $\gamma$ and $\partial Y$ are 
    in $\base(\mu)$, then we 
    have $\pi_{Z}(\mu)\subset B_{3}(\pi_{Z}(\partial Y))$ and 
    $\pi_{Y}(\mu)\subset B_{3}(\gamma)$. Thus $(\gamma,\pi_{Z}(\partial 
    Y))\subset B_{6}(\pi_{Y}(\mu) \times \pi_{Z}(\mu))$. Thus we have shown that 
    $\C(Y) \times \pi_{Z}(\partial Y) 
    \subset  B_{6}(\pi_{Y} \times \pi_{Z}(\M))$. The identical 
    argument with the 
    roles of $Y$ and $Z$ reversed completes the proof of the density 
    statement.    
\end{proof}

\begin{rmk} In the above proposition we see that except in the case
where $Y \cap Z =\emptyset$, there are uniform bounds for which the image of 
the map from $\C(S)$ to $\C(Y) \times \C(Z)$ lies in a neighborhood of 
the 
$\delta$--hyperbolic space formed from the complexes of curves for $Y$
and $Z$ ``joined together.''  The two complexes are glued together along a bounded 
diameter set in the
non-nesting case.  In the case of nesting $\C(Y)$ is glued to the link of a
point in $\C(Z)$; essentially this is done by taking a ``blow-up'' of 
$\C(Z)$ at the point $\partial Y\in\C(Z)$.  
$\delta$--hyperbolicity in the nesting case follows from a more 
delicate analysis, see \cite{Behrstock:thesis}.
\end{rmk}

\begin{figure}[ht!]
\labellist\small
\pinlabel {$\pi_{Z}(\partial Y)$} [bl] at 192 77
\pinlabel {$\pi_{Z}(\partial Y)$} [tl] at 260 73
\pinlabel {$\pi_{Y}(\partial Z)$} [br] at 170 38
\hair5pt
\pinlabel {$Z\cap Y =\emptyset$} [t] at 52 22
\pinlabel {$Z\cap Y\neq\emptyset$} [t] at 170 22
\pinlabel {non-nested} [t] <0pt,-10pt> at 170 22
\pinlabel {$Y\subset Z$} [t] at 287 22
\endlabellist
\centerline{\includegraphics[width=4.4truein]{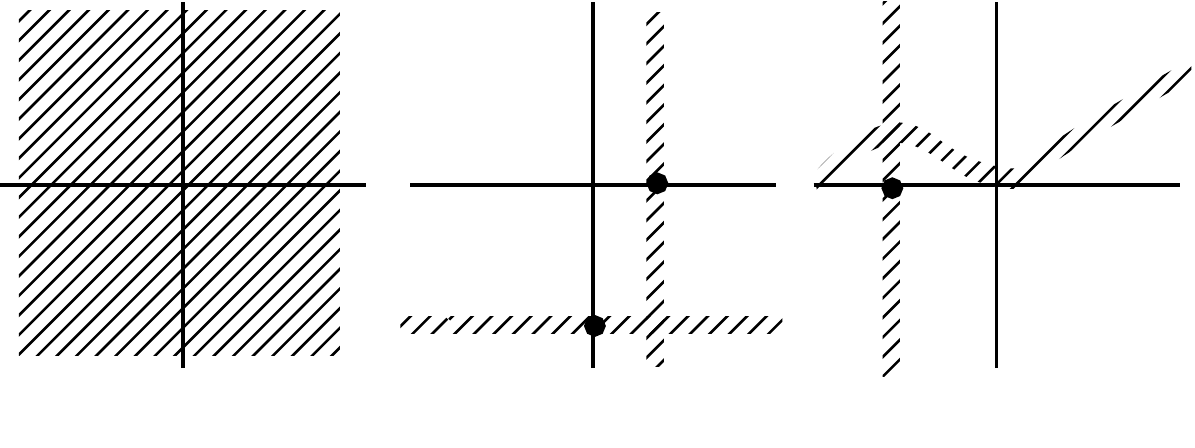}}
\caption{A cartoon of the
possible images of projections from complete markings into
$\C(Z)\times\C(Y)$}
\label{projpicts}
\end{figure}
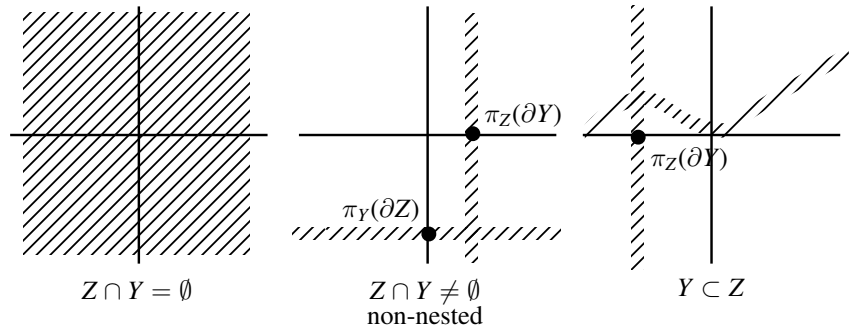

\fullref{geomprojest} (Projection estimates, geometric version) 
is summarized in \fullref{projpicts}. For a (slightly) more 
accurate picture of the behavior in the nesting case, see also 
\fullref{nestedprojpicts}.
 
\section{Hyperbolicity of the mapping class group and 
Teich\-m\"{u}ller space in the cases of low complexity}
    \label{sectionhyperbolicitylowcomplexity}

In this section we give new proofs of the following theorems.

\begin{thm}\label{pantshyperbolic}
    $\P(S_{1,2})$ and $\P(S_{0,5})$ are each $\delta$--hyperbolic.
\end{thm}

\begin{thm}\label{mcghyperbolic}
    $\M(S_{1,1})$ and $\M(S_{0,4})$ are each $\delta$--hyperbolic.
\end{thm}

\fullref{mcghyperbolic} is classical, since $\M(S_{1,1})$ is
isomorphic to $SL_{2}(\Z)$ and $\M(S_{0,4})$ to a finite subgroup of
$SL_{2}(\Z)$, hence both group are virtually free and thus easily seen to be 
$\delta$--hyperbolic. 
\fullref{pantshyperbolic} is a much deeper theorem, 
and was originally proved by Brock--Farb in \cite{BrockFarb:curvature} 
(see also \cite{Aramayona:thesis}). 
Although the Weil--Petersson metric has been known for some time to have 
negative curvature, it was only recently discovered in 
\cite{Huang:AsympFlat} that the 
sectional curvature is not bounded away from zero, even for the case of 
$\P(S_{1,2})$ and $\P(S_{0,5})$, thereby prohibiting a proof of 
\fullref{pantshyperbolic} by a comparison geometry 
argument.

We find especially tantalizing the phenomenon that these 
two results about two very different spaces can be proven 
simultaneously. This 
can be done since the property of low complexity which we use is 
that any two proper subsurfaces (with nontrivial curve complex) must 
overlap: a property of $S_{1,1}$ and $S_{0,4}$ which is also true 
for $S_{1,2}$ 
and $S_{0,5}$ when one adds the assumption that the subsurfaces are 
not annuli, a natural assumption when one is considering pants 
decompositions instead of markings (cf.\ \fullref{pantsalmostquasi}). 
Thus the proof of \fullref{pantshyperbolic} is 
obtained by rereading the proof for \fullref{mcghyperbolic} 
reading ``marking'' as ``pants decomposition,'' 
``hierarchy'' as ``hierarchy without annuli,'' 
and $\M(S)$ as $\P(S)$. The core arguments in this section will 
be revisited to prove \fullref{Fzertree}, a result which holds for 
surfaces of arbitrary complexity.

Our method of proof is to show that for every hierarchy path there 
exists an ``almost locally constant'' map sending the marking (or 
pants) complex to the hierarchy path. That this implies hyperbolicity is a 
consequence of a general result from \cite{MasurMinsky:complex1}; we 
also provide a new proof of this implication using asymptotic cones 
(\fullref{locallyconstantimplieshyperbolic}).

\subsection{A projection from markings to hierarchies}

Fix a surface $S$, a pair of complete clean markings $I$ and $T$, 
and a hierarchy $H$ connecting them.
By a hierarchy path, we mean a ``resolution'' of slices of the hierarchy into a 
sequence of markings separated by elementary moves; this is done by taking the
slices of $H$ and 
interpolating between them to get a path. The choice of a hierarchy path
is not canonical, but they are each $(K,C)$--quasi-geodesics with the 
quasi-isometry constants depending only on the topological type of the surface.
(See \cite{MasurMinsky:complex2}, 
especially Section 5 and the Efficiency of hierarchies Theorem of 
Section 6.)

We say that \emph{$Y$ appears as a large domain in $H$} if $Y$ is in 
the set
$$\G=\{Y\subseteq S: d_{Y}(I(H),T(H))>6M+ 4\delta'\}.$$
Where $M=\max\{M_{1}+2,M_{2}+3\}$ is 
the constant coming from \fullref{projest} (Projection estimates), 
note in 
particular that this constant is sufficiently large to satisfy the 
hypothesis of \fullref{MMtwLLL} (Large 
link \cite{MasurMinsky:complex2}) and 
$\delta'=4\delta+5$ where $\delta$ denotes the maximum of the
$\delta$--hyperbolicity 
constants for $\C(Y)$ where $Y\subseteq S$. Alternatively, when dealing 
with the pants complex we consider:
$$\G'=\{Y\in\G: Y\mbox{ is not an annulus}\}.$$
Letting $g_{H,Y}$ denote the geodesic segment supported on $Y$ in the 
hierarchy $H$, 
for each $Y\in \G$ we have a map
$\M(S)\xrightarrow{\pi_{Y}} \C(Y) \xrightarrow{r_{Y}} g_{H,Y}$ where $r_{Y}$ 
is the closest point(s) projection from $\C(Y)$ to $g_{H,Y}$.
We denote the composition
$$p_{Y}=r_{Y}\circ\pi_{Y} \co \M(S) \rightarrow g_{H,Y},$$ or
just $p$ when the surface $Y$ is understood. Notice that 
$\delta$--hyperbolicity of 
$\C(Y)$ (via \fullref{lemmaeasythintriangleargument}) 
combined with the fact that $\diam_{Y}(\pi_{Y}(\mu))\leq 3$ 
(\fullref{MMtwlpl}) 
imply that for any $\mu\in\M(S)$, the set 
$p_{Y}(\mu)$ has diameter at most $\delta'$. Accordingly, it follows
 that  maps $p_{Y}$ are coarsely distance decreasing 
in the sense that:
\begin{align}
    \label{eqcoarselydistancedecreasing} 
    d_{\C(Y)}(p_{Y}(\mu),p_{Y}(\nu))\leq  
    d_{\C(Y)}(\pi_{Y}(\mu),\pi_{Y}(\nu))+2\delta'.\\
\mathcal{L}(\mu)=\{Y\in \G: d_{Y}(p_{Y}(\mu),T(Y))<3M+2\delta'\}
\tag*{\hbox{Define }}\\
\mathcal{R}(\mu)=\{Y\in \G: d_{Y}(p_{Y}(\mu),I(Y))<3M+2\delta'\}. 
\tag*{\hbox{and}}
\end{align}
Usually the choice of marking $\mu$ will be fixed and we will drop 
the $\mu$ from the notation, just writing $\mathcal{L}$ and 
$\mathcal{R}$.

The next lemma summarizes some basic facts we will need about these 
sets.

\begin{lem}\label{lemmalowcomplexity}\begin{enumerate}
    \item $\mathcal{L}\cap\mathcal{R}=\emptyset$

    \item Let $Y,Z\in \G$. If $Y \prec_{t} Z$ and $Y$ and $Z$ overlap, 
    then either $Y\in \mathcal{L}$ or $Z\in \mathcal{R}$.
    \item If $S$ is $S_{1,1}$ or $S_{0,4}$
     then the set
     $\G\setminus (\mathcal{L}\cup\mathcal{R})$ consists of at most one 
     surface.
     \item If $S$ is $S_{1,2}$ or $S_{0,5}$ then 
     $\G'\setminus (\mathcal{L}\cup\mathcal{R})$
     consists of at most one surface.
\end{enumerate}
\end{lem}

\begin{proof}\ \

(1)\qua If $Y\in\mathcal{L}$ then 
    $d_{Y}(p(\mu),T)<3M$. If $Y$ is also in $\mathcal{R}$ then  
    $d_{Y}(p(\mu),I)<3M$. But together these imply that the total 
    distance $H$ travels through $Y$ is less than $6M$ 
    contradicting our assumption that $Y$ is a large domain.

(2)\qua Suppose not. Then we have $Y \prec_{t} Z$,  $Y\notin 
    \mathcal{L}$, and $Z\notin \mathcal{R}$.
    
    Since $Y$ and $Z$ are large domains with $Y \prec_{t} Z$, by 
    \fullref{MMtwOPL} 
    we have: $$d_{Y}(\partial Z,T(H))\leq M_{1}+2 
    \mbox{ and } d_{Z}(\partial Y,I(H))\leq M_{1}+2.$$
    $Y\notin \mathcal{L}$ implies that 
    $$d_{Y}(p(\mu),T(Y))\geq 3M+2\delta'.$$ Similarly,  $Z\notin \mathcal{R}$ implies 
    that $$d_{Z}(p(\mu),I(Z))\geq 3M+2\delta'.$$    
    Putting these facts together with the result from 
    \cite{MasurMinsky:complex2} that for any domain\break $d_{D}(I(D),I(H))\leq M_{1}$,
    yields the two inequalities: $$d_{Y}(p(\mu),\partial Z)>3M+2\delta'-M- 
    (M_{1}+2)>M+2\delta'$$
    $$d_{Z}(p(\mu),\partial Y)>3M+2\delta'-M-(M_{1}+2)>M+2\delta'.
    \leqno{\hbox{and}}$$
    The fact that $p_{Y}$ and $p_{Z}$ are coarsely distance 
    decreasing implies, via Equation~\ref{eqcoarselydistancedecreasing}, 
     that both $d_{Y}(\mu,\partial Z)>d_{Y}(p(\mu),\partial 
     Z)-2\delta'>M$ and 
    $d_{Z}(\mu,\partial Y)>d_{Z}(p(\mu),\partial Y)-2\delta'>M$.
    But this is impossible, since together these inequalities contradict
    \fullref{projest} (Projection estimates).
    
(3)\qua This follows from part 2, since in each of 
    these surfaces any two subsurfaces overlap and thus are time ordered 
    by Lemma~4.18 of \cite{MasurMinsky:complex2}.

(4)\qua Again this follows from part 2, since here any two 
    non-annular surfaces must overlap.
\end{proof}

Parts 3 and 4 of the above proposition suggests the following construction, 
which we give for the marking complexes for $S_{1,1}$ and $S_{0,4}$ 
and also for the pants complexes of $S_{1,2}$ and $S_{0,5}$
(aspects of this construction generalize to work on arbitrary surfaces, 
we discuss this further in \fullref{sectiontransversaltree}).
Given any pair of markings $I,T\in\M(S)$ and a hierarchy $H$ 
connecting them let $\widetilde{H}$ denote the 
set of markings associated to slices of $H$.
Below we construct a map $\phihat$ which maps elements of 
$\M(S)$ to (uniformly) bounded diameter subsets of $\widetilde{H}$ (in 
the cases we are considering,  $\widetilde{H}$ is metrized by a linear time 
ordering on slices).  $\widetilde{H}$  
be embedded into $\M(S)$ via a quasi-isometric embedding which can be 
extended by resolving the set of slices into a path (referred to as a 
\emph{hierarchy path}), the map $\phihat$ induces a map to 
(uniformly) bounded diameter subsets of any hierarchy path from $I$ to $T$.
The map $\phihat$ is defined to be the identity on markings in $\widetilde{H}$;
the following defines the map for  $\mu\in\M(S)\setminus \widetilde{H}$:

\begin{enumerate}
    \item If $\G\setminus (\mathcal{L}\cup\mathcal{R})=\{A\}\neq 
    \emptyset$, then define $\phihat(\mu)=(A,p_{A}(\mu))$.

    \item If $\G\setminus (\mathcal{L}\cup\mathcal{R})=\emptyset$, then 
    consider 
    $$\Lambda=\{v\in g_{H}: d_{\C(S)}(v,\mu)\leq d_{\C(S)}(w,\mu) \mbox{ 
    for any } w\in g_{H}\},$$
    ie, the set of points on the main geodesic $g_{H}$ which are closest 
    to $\mu$, namely $r_{S}(\mu)$. $\phihat(\mu)$ is defined to be the following sets of markings.
    \begin{enumerate}
        \item If any of the surfaces in $\mathcal{L}$ are component 
        domains of $g_{H}$ at a vertex in the set $\Lambda$, then 
        denote by $L$
	the rightmost (the last to appear with respect to the time ordering) 
	domain in
	$\mathcal{L}$, and define $\phihat(\mu)$ to contain each 
	$(L,v)$ where $v$ is any vertex of 
	$g_{H,L}$ within $3M+3\delta'$ of the terminal marking $T$.
    
        \item If a component domain corresponding to a vertex of 
        $g_{H}$ in the set $\Lambda$  is an element of $\mathcal{R}$ 
	then choose the leftmost element of this set (call this $R$). 
	Define $\phihat(\mu)$ to contain each of the markings $(R,v)$ 
	where $v$ is any vertex of $g_{H,R}$ within $3M+3\delta'$ of the initial
	marking $I$.
    
        \item For each domain $D$ which is a component domain of a 
        vertex in $\Lambda$ and which  
        is not in $\G$,  if $D$ supports a geodesic and it is 
        time ordered after the $L$ geodesic and before the $R$ geodesic 
        then  define $\phihat(\mu)$ to contain each $(D,v)$ 
	where $v$ is any vertex in $g_{H,D}$.
    \end{enumerate}
\end{enumerate}

Before proceeding, we mention  why $(A, p_{A}(\mu))$, 
$(L,v)$, and $(R,v)$, although not necessarily complete clean 
markings,  give rise to such 
(in the text we abuse notation and drop this correspondence from the 
discussion, thereby simply referring to  $(A, p_{A}(\mu))$, etc, as 
complete clean markings). In each case these consist 
of a proper subsurface $X$ and a point in $\C(X)$ (we 
unify the argument that these give  markings by writing $(X, p_{X}(\mu))$). 
For $S_{1,1}$ and 
$S_{0,4}$, a point of $\M(S)$ consist of a curve in $\C(S)$ and a 
transversal to that curve, so it is clear that our prescription above 
indeed describes a marking since the only nontrivial proper subsurfaces 
are annuli so $(X, p_{X}(\mu))$ refers to the marking $(\partial X, p_{X}(\mu))$. 
The case of $\P(S)$ for the surfaces 
$S_{1,2}$ and $S_{0,5}$ requires (only slightly) more justification. In 
these cases an element of $\P(S)$ is a pair of disjoint curves. Recall 
that the only subsurfaces of $S$ with nontrivial curve complex 
are once punctured tori and four punctured spheres (we ignore 
annuli when dealing with the pants complex).
Thus $p_{X}(\mu)\in\C(X)$
is also an element of $\C(S)$. Furthermore, $p_{X}(\mu)$ and $\partial 
X$ have distance 1 in $\C(S)$, so $(X, p_{X}(\mu))$ is taken to refer to the 
pants decomposition $(\partial X, p_{X}(\mu))$.

\begin{lem}\label{lemmadiameterbound} There is a uniform bound $D$ 
    (depending only on the  topological type of $S$), 
    so that for each $\mu\in\M$ the set $\phihat(\mu)$ has diameter 
    less than $D$.
\end{lem}

\begin{proof} First note that the set of points on the main 
    geodesic $g_{H}$ which are closest  
    to $\mu$ has diameter at most $\delta'$ by  
    \fullref{lemmaeasythintriangleargument}. 
    When $\G\setminus (\mathcal{L}\cup\mathcal{R})\neq\emptyset$ 
    there is only one such surface, as proven in 
    \fullref{lemmalowcomplexity}. In this case it follows from the 
    definition of $\Phi$ and 
    $\delta$--hyperbolicity of $\C(A)$ that $\Phi(\mu)$ is a subset of 
    $\widetilde{H}$ of diameter at most $\delta'$.

    Note that under $\Phi$ the marking $\mu$ can not
    project to anything time ordered before the rightmost element 
    $L\in\mathcal{L}$ as then \fullref{MMtwOPL}  would 
    force $p_{L}(\mu)$ to lie near the initial marking of $L$ and then we 
    would have either $L\notin \G$ or $L\in\mathcal{R}$, either way 
    contradicting \fullref{lemmalowcomplexity} which proves that 
    every surface in $\mathcal{L}$ is time ordered before every 
    surface of $\mathcal{R}$;  a similar argument gives the analogous 
    result for $\mathcal{R}$.

    So now $\Psi(\mu)$ consists of the 
    rightmost element of $\mathcal{L}$ which we call $L$, the leftmost element
    of $\mathcal{R}$ which we call 
    $R$, and all the rest of the small domains supporting geodesics 
    which are time ordered between $L$ and $R$ 
    (of which there are at most $\delta'$). 
    $\Phi$ was defined to be the union of the three set (any of which 
    are possibly empty):
    \begin{enumerate}
	\item     If $L\in\Lambda$ then $(L,v)$ where $v$ is any vertex of 
	$g_{H,L}$ within $3M+3\delta'$ of the terminal marking $T$.

	\item     If $R\in\Lambda$ then $(R,v)$ where $v$ is any vertex of 
	$g_{H,R}$ within $3M+3\delta'$ of the initial marking $I$.

	\item If $D\in\Lambda$, $D\notin\G$, and 
	$L\prec_{t}D\prec_{t}R$ then $(D,v)$ where $V$ is any vertex of 
	$g_{H,D}$.
    \end{enumerate}

    First, note that the diameter of the elements in item 1 is at most 
    $B(3M+3\delta')$, where $B$ is the 
    Lipschitz constant for the map from hierarchy slices to paths in 
    the marking complex. (Similarly for the markings in item 2.)
    Since $\C(S)$ 
    is $\delta$--hyperbolic  $\diam_{\C(S)}(\Lambda)<\delta'$,
    and thus if both $L$ and $R$ are in $\Lambda$ then 
    $d_{\C(S)}(L,R)< \delta'$. Also note that for any $D\in\Lambda$ as 
    described in item 3 we have that the diameter in $\M(S)$ of the 
    set of $(D,v)$ is less than $B(6M+4\delta')$ (since
    $6M+4\delta'$ is the threshold for being in $\G$). 
    Again, since the diameter in $\C(S)$ is bounded by 
    $\delta'$ we then see that the set of all markings in item three has 
    diameter in the  marking complex bounded by $\delta'\cdot 
    B(6M+4\delta')$.
    
    So now it follows that although $\Phi(\mu)$ consist of many slices: 
    $$\diam_{\M(S)}(\Phi(\mu)<2B(3M+3\delta'+ 
    \delta'(6M+4\delta')).\proved$$
\end{proof}

Note that since the map from slices to hierarchy paths is Lipschitz,  
the above lemma yields a map which sends points of $\M$  
to uniformly bounded diameter subsets of the hierarchy path, where the 
constants are independent of the hierarchy path.

In the next subsection we will show that these maps to hierarchy paths 
are ``almost locally constant'' off the hierarchy path, in the sense that 
there are uniform constants, so that large diameter balls around 
markings far from a hierarchy path are sent to small diameter subsets of 
the hierarchy path.

\subsection{Coarsely contracting projections}

We now describe the tool we use for proving $\delta$--hyperbolicity. 
It is a generalization to our context of Morse's Lemma on stability of 
quasi-geodesics in hyperbolic space.

\begin{defn}\label{defcontractionproperty}
    In a space $X$, we say a family of paths 
    $\mathcal{H}$ is \emph{transitive} if every pair of points in $X$ 
    can be connected by a path in $\mathcal{H}$. 
    We say $\mathcal{H}$ has the 
    \emph{coarsely contracting property}
    when $\mathcal{H}$ is a transitive family of paths in $X$ with the 
    property that  for every 
    $H\in\mathcal{H}$ there exists a map $\phihat_{H}\co X\to H$ 
    and constants $b$ and $c$ such that each for each  
    $\mu,\mu'\in X$ satisfying 
    $d_{X}(\mu,\mu')< b\cdot d_{X}(\mu, \phihat_{H}(\mu))$, then 
    $\diam(\phihat_{H}(\mu), \phihat_{H}(\mu'))<c$. 
    
\end{defn}

If considering $\M(S)$ then fix $S$ to be either $S_{1,1}$ or $S_{0,4}$, if 
considering $\P(S)$ then fix $S$ to be either $S_{1,2}$ or $S_{0,5}$ and use 
hierarchies-without-annuli.

\begin{lem}\label{lemmacontractionproperty}
    \noindent
    
    \begin{itemize}
	\item When  $S$ is either $S_{1,1}$ or $S_{0,4}$, then 
	the hierarchy paths form a coarsely contracting family of paths on $\M(S)$.

        \item When  $S$ is  either $S_{1,2}$ or $S_{0,5}$, then 
	the hierarchy-without-annuli paths 
	form a coarsely contracting family of paths on $\P(S)$.

    \end{itemize}
\end{lem}

\begin{proof} There exist a hierarchy connecting any pair of points in 
    $\M(S)$, which shows that  
the set of hierarchy paths form a transitive path family (similarly in 
$\P(S)$ for hierarchy-without-annuli paths).
Fix $H$ a hierarchy path between two points in $\M(S)$ and use $\phihat$ 
to denote $\phihat_{H}$.

We need to show: 
there exist constants $b$ and $c$ so that if 
$\mu,\mu'\in\M$ satisfy 
$d_{\M}(\mu,\mu')< b\cdot d_{\M}(\mu, \phihat(\mu))$, then 
$\diam(\phihat(\mu), \phihat(\mu'))<c$.

First notice that a key property of $\phihat$ is that there exists a 
constant $Q=3M+3\delta'$
so that every domain $A$ of a geodesic in $H$ (including the main surface $S$) 
has $\diam_{A}(\phihat(\mu)\cup p_{A}(\mu))<Q$; this fact is a step in 
the proof of \fullref{lemmadiameterbound}.

Fix a hierarchy $G$ from $\phihat(\mu)$ to $\mu$. 
Let $Z$ be a domain of $G$ of length at least $4Q$, an assumption 
which implies that:
\begin{equation}
    \label{eqZlarge} 
    d_{Z}(\mu, \phihat(\mu))>3Q+M.
\end{equation}
In order to compute 
$d_{\M}(\phihat(\mu),\phihat(\mu'))$ it is useful to calculate 
$d_{Y}(\mu, \mu')$ for each large domain $Y$ of $H$. Consider a
subsurface $Z$. 
If $Z$ is a large domain of $H$, then the definition of 
$\phihat$, \fullref{MMtwOPL}, and Equation~\ref{eqcoarselydistancedecreasing}  
together imply that  
$d_{Z}(\mu,\partial Y)\geq d_{Z}(\mu, \phihat(\mu))-5M-3\delta'$.
If $Z$ is not a large domain of $H$, 
then $d_{Z}(\partial Y, \phihat(\mu))<6M+4\delta'$; applying the triangle 
inequality then gives
$d_{Z}(\mu,\partial Y)\geq d_{Z}(\mu,\phihat(\mu))-\diam_{Z}(\phihat(\mu))-
d_{Z}(\phihat(\mu),\partial Y)-2\delta'$. 

Combining the two cases of the previous paragraph, for any surface $Z$, 
independent of whether or not it is a 
large domain in $H$, we have:
\begin{equation}
    \label{eqclosestpoint} 
    d_{Z}(\mu,\partial Y)\geq d_{Z}(\mu,\phihat(\mu)) -3Q.
\end{equation}
Combining the above  with  inequality~\ref{eqZlarge} we have   that 
$d_{Z}(\mu,\partial Y)> 3Q+M - 3Q=M$, which by the 
Projection estimates Theorem implies that $d_{Y}(\mu, \partial Z)<M$.

Now consider $d_{Z}(\mu', \partial Y)$: either $d_{Z}(\mu', \partial 
Y)>M$ for every large domain $Y$ of $H$ 
or  for some domain $Y$ of $H$ we have $d_{Z}(\mu', \partial Y)<M$.

In the first case, from \fullref{projest} (Projection estimates) 
we have 
$d_{Y}(\mu', \partial Z)<M$, which combines with the 
paragraph above to give $d_{Y}(\mu', \mu)<2M$. At which point it is easy to 
check that there is a uniform bound on 
$d_{\M}(\phihat(\mu),\phihat(\mu'))$---the bound comes since $d_{Y}(\mu', 
\mu)<2M$ implies that the pairs ($\mathcal{L}(\mu)$, $\mathcal{R}(\mu)$) 
and ($\mathcal{L}(\mu'),\mathcal{R}(\mu'))$ agree (except for a 
switch of at most one surface which then must have diameter not much 
larger than the threshold for $\G$) and that $\Lambda(\mu)$ 
is coarsely the same as $\Lambda(\mu')$. Carrying out this computation, 
one obtains $d_{\M(S)}(\phihat(\mu),\phihat(\mu'))<10M+6\delta'$.

In the other case, when there is some domain $Y$ of $H$ for which 
$d_{Z}(\mu', \partial Y)<M$, we can 
 combine this with equation~\ref{eqclosestpoint} to give:
\begin{equation}
    \label{eqestimatingmu} 
    d_{Z}(\mu,\mu')\geq d_{Z}(\mu, \phihat(\mu))-3Q-M.
\end{equation}
This tells 
us that whenever $d_{Z}(\phihat(\mu),\mu)>2(3Q+M)$, we have
\begin{equation}
    \label{eqdefinitemove} 
    d_{Z}(\mu,\mu')>\frac{1}{2}d_{Z}(\phihat(\mu),\mu)
\end{equation}
Thus, in domains $Z$ of $G$ for which
$d_{Z}(\phihat(\mu),\mu)$ is sufficiently large, either:\break
$d_{Z}(\mu', \partial Y)>M$,  
and thus $\phihat(\mu)$ and $\phihat(\mu')$ are close in the markings  complex, or 
 $d_{Z}(\mu,\mu')$ is a definite fraction of $d_{Z}(\phihat(\mu),\mu)$.

If we can show that a definite proportion of the distance travelled 
in $\M$ between $\phihat(\mu)$ and $\mu$, say $\frac{d_\M(\mu,\phihat(\mu))}{\alpha}$, 
 takes place in domains larger 
than a threshold $6Q+2M$ then we would be done by 
choosing the constant $b$ so that $b<\frac{d_\M(\mu,\phihat(\mu))}{\alpha}$. This would be enough,  as then $\mu'$ has 
no chance to move far enough away from $\mu$ to have $d_{Z}(\mu', \partial 
Y)<M$, where $Z$  is a large domain in $G$ and  $Y$ is a large domain in  $H$.

That there exists a uniform constant $\alpha>0$ for which the 
distance  $\frac{d_\M(\mu,\phihat(\mu))}{\alpha}$ 
occurs in large surfaces follows from the following counting argument.
The fact that the only domains of $S$ which support geodesics in 
$H$ are the surface $S$ and component domains of $g_{H}$ 
has as a consequence that if a hierarchy between $I$ and $T$ contains 
geodesics in $N+1$ different component domains, then $d_{\C(S)}(I,T)\geq N$. 
From this we see that either $\frac{1}{2}d_{\M(S)}(\phihat(\mu),\mu)$ 
occurs in large domains, or  the 
sum of the lengths of geodesics shorter than $6Q+2M$ is larger than 
$\frac{1}{2} d_{\M(S)}(\phihat(\mu),\mu)$. In the latter case, we then 
have $|g_{H}|>\frac{d_{\M(S)}(\phihat(\mu),\mu)}{2(6Q+2M)}$. Thus we 
always have that at least  $\frac{1}{2(6Q+2M)} d_{\M(S)}(\phihat(\mu),\mu)$ 
occurs in the large surfaces, as claimed.

Thus, choosing $b=\frac{1}{2(6Q+2M)}$ and  $c>10M+6\delta'$, 
we have proved the theorem.
\end{proof}

\subsection{Hyperbolicity}

In this section we prove Theorems~\ref{mcghyperbolic} and 
\ref{pantshyperbolic} using the contracting properties of hierarchy 
paths provided by \fullref{lemmacontractionproperty}.

In \cite{MasurMinsky:complex1} a contraction property is defined similar to 
 \fullref{defcontractionproperty} and is used to prove 
 $\delta$--hyperbolicity of $\C(S)$; that property motivated 
 the definition of the previous section. In this section, 
 using results about $\R$--trees we give a new proof that  such 
 contraction properties 
imply hyperbolicity, using techniques different than those of 
\cite{MasurMinsky:complex1}. This proof provides a new perspective
concerning the
relationship between stability of quasigeodesics and hyperbolicity.

\begin{thm}\label{locallyconstantimplieshyperbolic}
    If $X$ is a geodesic space with a family of $(K,C)$--quasi-geodesic 
paths $\mathcal{H}$ which 
have the coarsely contracting property, then $X$ is $\delta$--hyperbolic.
\end{thm}

Theorems~\ref{mcghyperbolic} and \ref{pantshyperbolic} follow 
 from this theorem since \fullref{lemmacontractionproperty} 
proves that in the low complexity cases $\M(S)$ (and $\P(S)$) 
the hierarchy paths form a coarsely contracting family of paths 
(respectively hierarchy-without-annuli paths), 
and \cite[Efficiency of hierarchies Theorem]{MasurMinsky:complex2} 
implies that hierarchy paths 
are $(K,C)$--quasi-geodesics, with the constants depending only on the 
topological type of the surface.

\begin{proof}
The method of proof is to show that in the asymptotic cone $\phihat$ 
induces a locally constant map to a certain family of bi-Lipschitz 
paths,  we then show that this implies that $\cone(X,\<d_{i}\>)$ 
is an $\R$--tree.
In particular we will show for each 
    pair of points $x,y \in\cone(X,\<d_{i}\>)$ there exists a 
    map $\Phi\co\cone(X,\<d_{i}\>) \to [x,y]$ such that:
\begin{enumerate}
    \item $\Phi|_{[x,y]}$ is the identity

    \item $\Phi|_{\cone(X,\<d_{i}\>)\setminus [x,y]}$ is locally constant.
\end{enumerate}

  By our hypothesis, for each $H\in\mathcal{H}$ there exists a map 
  $\phihat_{H}\co X\to X$ such that:
\begin{enumerate}
    \item For each $\mu\in X$ we have $\phihat_{H}(\mu)\in H$.

    \item There exist constants $b,c>0$ such that for any sequence of
    paths $H_{n}\in\mathcal{H}$ and points $\mu_{n}$ for which $r_{n}=
    d_{X}(\mu_{n},\pi_{H}(\mu_{n}))$ grows linearly, then\break
     $\diam_{X}(\phihat_{H_{n}}(B_{b\cdot r_{n}}(\mu_{n})))\leq c.$
\end{enumerate}

Fix three points 
$\bar{I}=\<I_{i}\>_{i\in\N}, \bar{T}=\<T_{i}\>_{i\in\N}, 
\bar{\mu}=\<\mu_{i}\>_{i\in\N}\in\cone(X,\<d_{i}\>)$.
Let $\bar{H}$ denote the path connecting $\bar{I}$ to 
$\bar{T}$ given by a rescaled sequence of paths $H_{i}$ 
connecting $I_{i}$ to $T_{i}$ and assume $\bar{\mu}\notin \bar{H}$. Since  
$\mathcal{H}$ has the coarsely contracting property, by definition we 
know that  
when $d_{X}(\mu_{i}, H_{i})$ grows linearly in $d_{i}$ 
(say it grows roughly as 
$g\cdot d_{i}$ for a constant $g$)
the ball with radius $b\cdot d_{X}(\mu_{i}, H_{i})$
around $\mu_{i}$ maps to a 
ball of radius at most $c$ in $H_{i}$. Thus in the rescaled metric 
space, $\frac{1}{d_{i}}\cdot X$, we have 
$d_{\frac{1}{d_{i}}\cdot X}(\mu_{i}, H_{i})$ is approximately $g$ 
and in this metric the ball $B_{b\cdot g d_{i}}(\mu_{i})$  maps to 
a set of diameter $\frac{c}{d_{i}}$ in $H_{i}$. 
Thus $\diam_{\cone(X,\<d_{i}\>)}(\<\phihat(B_{b\cdot g}(\mu_{i}))\>)=\ulimi 
\frac{c}{d_{i}}=0$, thereby showing that the map is locally 
constant.

This proves that any embedded path connecting  $\bar{I}$ and $\bar{T}$ must 
be a subset of $\bar{H}$. Since by hypothesis each $H\in\mathcal{H}$ is a 
$(K,C)$--quasi-geodesics, we have that $\bar{H}$ is the 
ultralimit as $i\to \infty$ of $(K,\frac{C}{d_{i}})$--quasi-geodesics  and is 
thus a $K$--bi-Lipschitz embedded path. In particular it is 
homeomorphically embedded, so the only subpath connecting its 
endpoints is the path itself---thereby proving  that there is a unique 
topological arc between any pair of points in $\cone(X)$.
It is worth remarking that
since $\cone(X)$ is a geodesics space, the path $\bar{H}$ which 
{\it a fortiori} need only be bi-Lipschitz embedded is actually 
the image of the geodesic between $\bar{I}$ and $\bar{T}$.

Mayer and Oversteegen's topological characterization of 
$\R$--trees \cite{MayerOversteegen:Rtrees} implies that in order to 
prove that $\cone(X)$ is an $\R$--tree it suffices to show that it is 
uniquely arcwise connected and locally path connected. The first 
property was shown above, the second is true for any geodesic space.

Thus each $\cone(X)$ is an $\R$--tree and thus \fullref{treeimplieshyper} 
implies that $X$ is $\delta$--hyperbolic.
\end{proof}

\section{Asymptotic cones of $\MCG$}\label{sectiontransversaltree}

When $\xi(S)\geq 2$, the mapping class group of $S$ is not 
hyperbolic and thus the results of the previous section clearly do 
not hold for arbitrary surfaces (similarly for the Teichm\"{u}ller 
space of $S$ when $\xi(S)\geq 3$). Nonetheless, in this section we 
show that negative curvature phenomena is present in certain 
directions. We begin by building a projection to mapping class groups 
of subsurfaces which we use to construct a subset which we call 
the strongly bounded subset, $F_{0}$, of the asymptotic cone. 
\fullref{Fzertree} is the main result of this section, it describes 
various topological properties of the subset $F_{0}$.
The proof of \fullref{Fzertree} is similar in spirit 
to the proof of Theorems~\ref{pantshyperbolic} 
and~\ref{mcghyperbolic}. The main difference, is that due to the lack 
of $\delta$--hyperbolicity in the spaces being considered in this 
section, we have no hope of showing that there is a 
coarsely locally constant projection to every hierarchy path. 
We consider a subcollection of hierarchy paths for which we build 
such projections; it the ultralimit of this collection of 
path that yield the set $F_{0}$.

\bigskip

Given any subsurface $Z\subset S$, we define a projection 
$\pi_{\M(Z)}\co\M(S)\to 2^{\M(Z)}$, which sends elements of $\M(S)$
to  subsets of $\M(Z)$. Given any $\mu\in\M(S)$ we build  
a marking on $Z$ in the following way. Choose an 
element $\gamma_{1}\in\pi_{Z}(\mu)$, and then recursively choose 
$\gamma_{n}$ from $\pi_{Z\setminus \cup_{i<n}\gamma_{i}}(\mu)$, 
for each $n\leq \xi(Z)$. Now take these $\gamma_{i}$ to be 
the base curves of a marking on $Z$. For each 
$\gamma_{i}$ we define its transversal $t(\gamma_{i})$ to 
be $\pi_{\gamma_{i}}(\mu)$. This process yields a complete marking 
on $Z$, since there are $\xi(Z)$ base curves and each has a
transversal; work of Masur and Minsky \cite[Lemma
2.4]{MasurMinsky:complex2}  
shows that any complete marking uniquely defines a finite 
subset of complete clean markings of diameter at most 6.
Thus this construction yields a bounded diameter subset of 
$\M(Z)$, but arbitrary choices in choosing the $\gamma_{i}$ 
were made along the way; we define 
$\pi_{\M(Z)}(\mu)$ to be the union of all possible markings built 
following this process, in the following lemma we show that 
this set has bounded diameter.

\begin{lem}\label{diamboundmarkingproj}
    There exist a uniform bound depending only on $S$, which bounds
the diameter of 
    $\diam_{\M(Z)}(\pi_{\M(Z)}(\mu))$ for any $\mu\in\M(S)$ and
$Z\subset S$.
\end{lem}

\begin{proof} Fix $\mu\in\M(S)$, $Z\subset S$, and
$\nu,\nu'\in\pi_{\M(Z)}(\mu)$. 
    It suffices to consider $\nu$ with base curves,
    $\gamma_{i}$, and transversals, $\pi_{\gamma_{i}}(\mu)$,  as
given above since 
    \cite[Lemma 2.4]{MasurMinsky:complex2} 
    implies that each such complete marking 
    has a bounded diameter subset of $\M(Z)$ associated with
    it (similarly for $\nu'$). We consider 
the projections of $\nu$ into  $\C(Y)$ for each  $Y\subseteq 
Z$ and show that the image under these projections are close to
$\pi_{Y}(\mu)$. 
We proceed inductively, by considering the smallest $i$ for which 
$\pi_{Y}(\gamma_{i})\neq\emptyset$.

If $i=1$, then since 
$d_{\C(Z)}(\gamma_{1},\mu)\leq 1$, \fullref{MMtwlpl} 
(Lipschitz projection) implies that 
$d_{\C(Y)}(\gamma_{1},\mu)\leq 3$ 
for any $Y\subset Z$. Thus we obtain $d_{\C(Y)}(\nu,\mu)\leq 3$. 

We now assume the inductive hypothesis that for a fixed $n$ we have 
$d_{\C(Y)}(\nu,\mu)\leq 3$ for all $Y\subset Z$ with 
$\pi_{Y}(\gamma_{i})\neq\emptyset$ for some $i<n$. Our goal is to
show 
that if $n$ is the smallest index of a curve
$\gamma_{n}\in\base(\nu)$ 
with $\pi_{Y}(\gamma_{n})\neq\emptyset$, then
$d_{\C(Y)}(\nu,\mu)\leq 3$. Since $\gamma_{i}\cap 
Y=\emptyset$ when $i<n$, we have 
$d_{\C'(Y)}(\pi'_{Y}(\mu),\pi'_{Y}(\gamma_{n}))\leq 1$. Thus by 
\fullref{MMtwlpl} we have 
$d_{\C(Y)}(\pi_{Y}(\mu),\pi_{Y}(\gamma_{n}))\leq 3$, and hence 
$d_{Y}(\mu,\nu)\leq 3$. Thus for any surface $Y\subset Z$ with 
$\pi_{Y}(\base(\nu))\neq\emptyset$ we have $d_{Y}(\mu,\nu)\leq 3$.

It remains 
to deal with the case where $Y$ is an annulus with base curve 
$\gamma\in\base(\nu)$, but here it follows from the definition of 
$\pi_{\M(Z)}(\mu)$ that $\pi_{Y}(\mu)=\pi_{Y}(\nu)$, thereby proving 
that for any $\nu\in \pi_{\M(Z)}(\mu)$ and any $Y\subset Z$ we 
have $d_{Y}(\mu,\nu)\leq 3$. This implies that 
and any $Y\subset Z$ we have  $\diam_{Y}(\nu,\nu')\leq 9$. 

Now, \fullref{mcgalmostquasi} shows 
that $d_{\M(Z)}(\nu,\nu')\sim \sum_{\substack{Y\subseteq Z \\ 
d_{Y}(\mu,\nu)>t}}d_{Y}(\nu,\nu')$, but we just showed that when the 
threshold is taken larger than 9 then  
the sum on the right side is zero. Thus if we choose a threshold
$t>9$ 
in \fullref{mcgalmostquasi}, then $d_{\M(Z)}(\nu,\nu')$ is
bounded 
by the additive quasi-isometry constant in the conclusion of that
theorem.
\end{proof}

As in the case of the maps $\pi_{Z}$, we abbreviate
$d_{\M(Z)}(\pi_{\M(Z)}(\mu_{i}),\pi_{\M(Z)}(\nu_{i}))$ by writing
$d_{\M(Z)}(\mu_{i},\nu_{i})$.

\begin{defn}\label{defFze}
Define a pair of points $(\mu,\nu)$  in $\coneM(S)\times \coneM(S)$ 
to be of \emph{bounded growth} if there exists a 
constant $D$ for which: 
$$\ulimi\sup\{d_{\M(Z)}(\mu_{i},\nu_{i}):Z\subsetneq S\}\leq D<\infty.$$
For our fixed basepoint $0\in\M(S)$, let $F_{0}$ denote the set of 
$\nu$ such that the pair $(0,\nu)$ is of bounded growth.  
We use the term \emph{strongly bounded} to refer to the elements 
of $F_{0}$.
\end{defn}

\begin{rmk} 
    Fixing a sequence of scaling factors $d_{n}$ and an ultrafilter  
    $\omega$, one can consider 
    those sequences $D(n)$ for which $\ulim \frac{D(n)}{d_{n}}=0$: 
    such sequences are said to have \emph{sublinear growth (with 
    respect to the scaling constants $d_{n}$ and ultrafilter $\omega$)}. 
    This definition allows one to consider, for instance, the set of points 
    $\mu\in\coneM(S)$ for which the sequence 
    $\{d_{\M(Z)}(\mu_{i},0):Z\subsetneq S\}$ grows 
    sublinearly, we call these points \emph{strongly sublinear 
    (for the choices $d_{n}$ and $\omega$)}.
    (A special case of this is for the  scaling constants 
    $d_{n}=\frac{1}{n}$, indeed this class of asymptotic cones was 
    used in many of the original papers in the subject, eg, 
    \cite{Gromov:PolynomialGrowth,Gromov:Asymptotic,DriesWilkie,Kapovich:GGTbook}.)

    We note that the remaining 
    results and proofs of this section hold mutatis mutandis replacing 
    ``strongly bounded'' by ``strongly sublinear.'' This is often 
    useful, but in the sequel we use the  
    strongly bounded hypothesis in order to have a condition which 
    holds in all asymptotic cones simultaneously. 
\end{rmk}

\begin{exa}\label{PAsublinear} If $\rho\in\MCG(S)$ is pseudo-Anosov, 
    then the sequence $\<\rho^{n}\>_{n}$ has distance from the 
    origin, $0$,  growing at a linear rate. At the same time, 
    there is a uniform bound $M$ so that for all $n\in\N$ and 
    $Y\subseteq S$ one has $d_{Y}(\rho^{n}, 0)\leq M$. 
    This shows that $|F_{0}|>1$. This is an important example, as in 
    this case, up to some quasi-isometry constants, 
    all the mapping class group distance is accounted for by 
    distance in the complex of curves of the surface $S$.
\end{exa}

This notion of degrees of boundedness  
provides an interesting stratification of $\coneM$;
Theorems~\ref{Fzertree} and \fullref{cutpointstw} are examples of 
how this stratification provides information about the mapping class
group. These results are refined further in Behrstock--Minsky's proof 
of the Rank Conjecture \cite{BehrstockMinsky:rankconj}.

In \fullref{sectionhyperbolicitylowcomplexity} we used
projections 
from $\M(S)$ to hierarchy paths and showed that these maps had a
strong 
contraction property. In fact, these maps were sufficiently 
contracting so  that for each hierarchy path, in the asymptotic cone
we 
obtained a  locally constant map from  $\coneM(S)$ to a bi-Lipschitz
path. 
From this we concluded that in the low complexity cases 
$\coneM(S)$ is an $\R$--tree. This argument relied heavily on the fact 
that these surfaces had low topological complexity, nonetheless we
extend the essential parts of that argument to produce strongly 
contracting maps in mapping class groups of surfaces of arbitrary
complexity. In so doing, the following result shows that 
bi-Lipschitz flats in $\coneM$ intersects $F_{0}$ in a very simple way. 
This answers and extends Conjecture~6.4.6 of \cite{Behrstock:thesis}. 
This contraction property also shows that given a pseudo-Anosov 
element $\phi$ and sufficiently large $n\in\N$, then any path in 
$\M(S)\setminus B_{|\phi^{2n}|}(1)$ has length which is at least 
quadratic in $|\phi^{2n}|$. This holds since the contraction property 
shows that linear length paths in this space move at most a uniformly 
bounded amount when projected to the quasi-geodesic from 
$\phi^{-n}(1)$ to $\phi^{n}(1)$, therefore any path in 
$\M(S)\setminus B_{|\phi^{2n}|}(1)$ connecting 
$\phi^{-2n}(1)$ to $\phi^{2n}(1)$ must contain at least linearly many 
disjoint linear length subpaths. 
Except when it is virtually free, $\M(S)$ is
thick of order~1 \cite{BehrstockDrutuMosher:thick} and thus has
at most quadratic divergence. Hence, except 
in those cases where it is virtually free,
the mapping class group has a quadratic divergence function. 
See the next subsection for some further topological
corollaries and a discussion of other applications.

 \begin{thm}\label{Fzertree} $F_{0}$ satisfies the following
properties:
     \begin{enumerate}
	 \item\label{Fzeparton} $F_{0}$ is geodesically convex.

	 \item\label{Fzeparttw} Any distinct pair of points in
	 $F_{0}$ are connected by a unique path in
	 $\coneM$. In particular,  $F_{0}$ is an $\R$--tree.

	 \item \label{Fzepartth} In $\coneM$, any bi-Lipschitz flat of 
	 rank larger than one intersects $F_{0}$ in at most
	 one point.
     \end{enumerate}
 \end{thm}

 Before proving this theorem we establish some notation and
 provide a lemma which we shall need in the proof.
 
 Similar to the definition in 
 \fullref{sectionhyperbolicitylowcomplexity}, for a 
 hierarchy $G$  we define the \emph{large domains} of $G$ as:
 $$\G_{G,C}=\{Y\subseteq S: d_{Y}(I(G),T(G))> C \geq 3M\}.$$
We will drop the $C$ from the notation when the choice of $C$ 
is irrelevant.
 Also, we define the \emph{initial large domains} by: 
 $$I_{G}=\{Y\in\G_{G}: \mbox{Each }Y'\in\G\mbox{ which overlaps }
 Y\mbox{ has } Y \prec_{t} Y' \}.$$ From finiteness of $\G$ it follows
 $I_{G}\neq\emptyset$; since any two domains of $I_{G}$ are disjoint 
 or nested there is a
 bound on $|I_{G}|$ depending only on the topological type of~$S$.

  \begin{lem}\label{lemmalargesetsexists}
      If $Y'\in\G_{G}\setminus I_{G}$ then there exists $Y\in I_{G}$
such that $Y'$  overlaps~$Y$.
  \end{lem}

  \begin{proof}
      Letting $Y'\notin I_{G}$, then there must exist a domain
$Y_{1}\in\G_{G}$ which
      overlaps $Y'$ and for which $Y_{1}\prec_{t}Y'$. Either
$Y_{1}\in I_{G}$
      in which case we are finished, or we can find a surface
      $Y_{2}\in\G_{G}$ which overlaps $Y_{1}$ and for which
      $Y_{2}\prec_{t}Y_{1}$. Since $Y_{1}\in \G$,
      \fullref{MMtwOPL} implies that $d_{Y_{1}}(\partial Y',
\partial
      Y_{2})>M$ and thus $Y'$ must overlap $Y_{2}$. Repeating this
      process finitely many times (since $\G_{G}$ is finite) yields a
      sequence $Y \prec_{t}\ldots\prec_{t} Y_{2}\prec_{t}Y_{1}
      \prec_{t} Y'$ with $Y\in I_{G}$ and $Y$ overlapping $Y'$.
  \end{proof}

  We now proceed to the proof of \fullref{Fzertree}.
 
 \begin{proof}  Fix two points $x=\<x_{i}\>,y=\<y_{i}\>\in F_{0}$.
     Fix a hierarchy $H_{i}$ connecting $x_{i}$
     to $y_{i}$ and let $[x_{i},y_{i}]$ denote the corresponding 
     hierarchy path between them.  We write $[x,y]$  to denote the 
     bi-Lipschitz path in $\coneM(S)$
     obtained by taking an ultralimit of these quasi-geodesic paths.

(1)\qua We will now show that every point of  $[x,y]$ is in $F_{0}$.
     By the definition of $F_{0}$, for every
     $\<\alpha_{i}\>\in \seq\setminus \<S\>$ we have
     $$\ulim (\max\{ 
     d_{\M(\alpha_{i})}(x_{i},0),
     d_{\M(\alpha_{i})}(y_{i},0)\})<\infty.$$
     Thus, in particular, 
     for each $\<\alpha_{i}\>\in\seq\setminus \<S\>$  we have 
     \begin{equation} \label{eqfzeparton}
	 \ulim\frac{1}{d_{i}}d_{\M(\alpha_{i})}(x_{i},y_{i})=0.
     \end{equation}
          We note the following useful property of hierarchy paths:
     if $[\mu,\nu]$ is a hierarchy path from $\mu$ 
     to $\nu$ and $\rho\in [\mu,\nu]$, then for all $Y\subsetneq S$
one has 
     $d_{\C(Y)}(\nu,\rho)\leq d_{\C(Y)}(\nu,\mu)+3M$ (where 
     $M$ is the constant of \fullref{projest},  
     depending only on $S$).
     Invoking \fullref{mcgalmostquasi} and 
     \fullref{diamboundmarkingproj} this implies that 
     $d_{\M(Y)}(\nu,\rho)< Kd_{\M(Y)}(\nu,\mu) +C$, for 
constants $K$ and
     $C$ which depend only on $S$.

     By definition of the asymptotic cone, if  $\<w_{i}\>\in[x,y]$, 
     then for  an $\omega$--large subsequence of indices the distances 
     $d_{\M(S)}(w_{i},[x_{i},y_{i}])$ grow slower than the sequence
$d_{n}$; 
     as a consequence of the 
     previous paragraph we see that for any  $\<w_{i}\>\in[x,y]$ 
     and any  $\<\alpha_{i}\>\in\seq\setminus \<S\>$ we have 
     $d_{\M(\alpha_{i})}(w_{i},y_{i})\leq
     Kd_{\M(\alpha_{i})}(x_{i},y_{i})+C$
     for $\omega$--almost every $i$.
     Together with equation~\ref{eqfzeparton} this implies that for 
     any $\<\alpha_{i}\>\in\seq\setminus \<S\>$ we have 
     $\ulimif d_{\M(\alpha_{i})}(y_{i},w_{i})=0$.
     Namely, we have shown that 
     $[x,y]\subset F_{0}$, ie, $F_{0}$ is geodesically convex.

(2)\qua This argument follows the general outline used to prove 
     \fullref{lemmacontractionproperty}.
     Consider the
     family $\mathcal{H}$  consisting of ultralimits of  rescaled
     hierarchy paths connecting pairs
     of points in $F_{0}$. We shall show that for
     any $x,y\in F_{0}$ and any  path $[x,y]\in \mathcal{H}$ 
     there is a continuous map $\Phi\co\coneM(S)\to[x,y]$ such that:
     \begin{enumerate}
	 \item $\Phi|_{[x,y]}$ is the identity
	 \item $\Phi$ is locally constant outside $[x,y]$.
     \end{enumerate}
     Once this is established
     \fullref{locallyconstantimplieshyperbolic}, 
     which we previously used in conjunction with
\fullref{lemmacontractionproperty} 
     to prove \fullref{pantshyperbolic}, implies that $F_{0}$ is
an $\R$--tree.
     The two key differences between this argument and that of
     \fullref{lemmacontractionproperty} are in
     the way in which the map $\Phi$ is constructed and then the
proof 
     that this map is locally constant---since the present
construction
     of $\Phi$ works for surfaces of
     arbitrary complexity we must rely on the property of being
     strongly bounded for both 
     these steps. 
     Similar to the construction of $\Phi$ in
     \fullref{sectionhyperbolicitylowcomplexity}, 
     we define our projection by considering a sequence of
hierarchies and maps to subsets of
     their slices and then perform a limiting procedure.

     \subsubsection*{$\Phi$ is well-defined}$\phantom{9}$
     Fix $\mu\in\coneM$ and let $g_{H_{i}}\subset\C(S)$ denote
     the main geodesic of $H_{i}$.
     By $\delta$--hyperbolicity of $\C(S)$, {\it a la}
     \fullref{lemmaeasythintriangleargument}, the closest point 
     projection of $\pi_{S}(\mu_{i})$ to $g_{H_{i}}$ 
     determines a subset $C_{\mu_{i}}\subset g_{H_{i}}$ whose
     diameter is uniformly bounded by a constant $K(\delta)$, which
     depends only on the $\delta$--hyperbolicity constant of $\C(S)$.
     Define  $\phihat(\mu_{i})$ to be the set of all marking in
$[x_{i},y_{i}]$
     whose base curves  include any element $\gamma\in C_{\mu_{i}}$.

     Note that given any element $\gamma\in C_{\mu_{i}}$, 
     by strong boundedness there are at most a bounded number 
     of elements of  
     $\C(\gamma)$ which appear as transverse curves to $\gamma$ in
     some marking of $H_{i}$.
     Further, strong boundedness also implies there is a uniform
bound on $\diam_{\C(S\setminus \gamma)}([x_i,y_i])$. Thus,
inductively building the subset      $\Phi(\mu)\subset[x_i,y_i]$ one
sees that there is a uniform bound, $K$, on its diameter.
   
   In particular, since for each $\mu_i\in[x_i,y_i]$ one has 
   $\mu_i\in\Phi(\mu_i)$, the above uniform bound implies that  
       $\diam_{\coneM}(\<\Phi(\mu_i)\>,\<\mu_i\>)=\ulimif(K)=0$. 
       Thus one obtains that $\Phi(\mu)$ is always a point, this also 
       shows that $\Phi$   restricts to the identity map on  $[x,y]$.

     \subsubsection*{$\Phi$ is locally constant off $[x,y]\phantom{9}$}
     Fix a point $\mu\in \coneM(S)\setminus [x,y]$ and a 
     point $\mu'\in \coneM(S)$ with $d(\mu,\mu')<c\cdot
     d(\mu,\Phi(\mu))$ for a positive constant $c$, depending only on
     the topological type of $S$, to be determined by the
     following argument.
     We want to show that for any such pair one obtains
$\Phi(\mu)=\Phi(\mu')$.

     Fix a representative for $\mu=\<\mu_{i}\>$ and for each $i$
consider
     a hierarchy $G_{i}$ whose initial marking is a point of
$\phihat(\mu_{i})$ 
     and with terminal marking $\mu_{i}$---this
     involves an arbitrary choice of element of $\phihat(\mu_{i})$,
     but the fact that $\diam (\phihat(\mu_{i}))$ is uniformly 
     bounded implies that 
     any choices here yield the same point in the asymptotic cone. As
     usual, $g_{G_{i}}$ will denote the main geodesic of $G_{i}$.
     
     The remained of the argument will be divided into two cases
     depending on the length of  $g_{G_{i}}$. 
     Let $K$ denote a fixed constant coming from
     \fullref{mcgalmostquasi} ($K$ is chosen by taking the
     quasi-isometry constant given by \fullref{mcgalmostquasi}
     corresponding to a 
     threshold $t>4M+D$, where $M$ is the constant given by
     \fullref{projest} and $D$ is the constant quantifying the
     strong boundedness of the pair $(x,y)$, as given in
     \fullref{defFze}).
     
     \subsubsection*{Case (i)}
     $|g_{G}|>\frac{1}{K}d(\mu,\Phi(\mu))$\quad
     More precisely, here we assume that for $\omega$--almost
     every $i$ the length
     of $g_{G_{i}}$ is larger than $\frac{1}{K}
     d_{\M(S)}(\mu_{i},\phihat(\mu_{i}))$.

     Let $c<\frac{1}{2K}$. Hyperbolicity of $\C(S)$, then implies 
     that there is a uniform bound on 
     $\diam_{\C(S)}(C(\mu_i)\cup C(\mu'_i))$ 
     (see \fullref{curvecomplexhyperbolic} and 
     \fullref{lemmaeasythintriangleargument}).  The proof that 
     $\Phi$ is well defined, then implies that  $\Phi(\mu)=\Phi(\mu')$.

     \subsubsection*{Case (ii)}
      $|g_{G}|\leq\frac{1}{K}d(\mu,\Phi(\mu))$\quad We now consider the 
     cases not covered by Case~(i), namely we assume that for
     $\omega$--almost every $i$ the length of $g_{G_{i}}$ is at most\break
     $\frac{1}{K} d_{\M(S)}(\mu_{i},\phihat(\mu_{i}))$.  The choice of 
     $K$ implies
     that for $\omega$--almost every $i$, any hierarchy between
     $\mu_{i}$ and $\phihat(\mu_{i})$ must have a domain $Y\subsetneq
     S$ for which $d_{Y}(\mu_{i},\phihat(\mu_{i}))>4M+D$, where $M$ is
     the constant given by \fullref{projest} and $D$ is the 
     strongly boundedness constant for the pair $(x,y)$.  Fix a constant
     $0<c<1$ and let $\mu_{i}=\mu_{i,0},
     \mu_{i,1},\mu_{i,2},\ldots,\mu_{i,k}=\mu_{i}'$ be an elementary move
     sequence of a hierarchy path in $\M(S)$ with 
     $k< c d_{\M(S)}(\mu_i,\phihat(\mu_i))$.

     For any surface $Y\subsetneq S$, define $$J_{i,Y}=\{j\in[0,k]:
     d_{Y}(\mu_{i,j},\mu_{i})>4M+D\mbox{ and }
     d_{Y}(\mu_{i,j},\phihat(\mu_{i}))>4M+D\}.$$
     We now give a brief digression to study properties of $J_{i,Y}$. 
     Suppose $j\in J_{i,Y}$, $l\in J_{i,Z}$, and $Y,Z\subsetneq S$
     overlap.  Then, $Y$ and $Z$ are both large domains for $G_{i}$, and
     in the time ordering on the hierarchy $G_{i}$ either $Y\prec_{t}Z$ or
     $Z\prec_{t}Y$, for the remainder of the paragraph we assume the
     first ordering (in the other case an analogous argument can
     easily be provided by the reader).  \fullref{MMtwOPL} shows
     that $Y\prec_{t}Z$ implies $d_{Y}(\partial Z, \mu_{i})<M$ and
     $d_{Z}(\partial Y, \phihat(\mu_{i}))<M$; the first inequality
     combines with the definition of $J_{i,Y}$ to imply that
     $d_{Y}(\mu_{i,j},\partial Z)>3M+D$.  \fullref{projest} then
     implies that $d_{Z}(\mu_{i,j},\partial Y)< M$, which 
     combines with the second inequality of the
     previous sentence to show that $d_{Z}(\mu_{i,j},
     \phihat(\mu_{i}))<2M$.  Thus $j\notin J_{i,Z}$ and we have shown
     that for any $Y,Z\subsetneq S$ we have $J_{i,Y}\cap
     J_{i,Z}=\emptyset$.  Thus $J_{i,Y}$ and $J_{i,Z}$ share a point
     only if they are disjoint or nested, a straightforward
     topological count gives a bound of $2\xi(S)$ for the number of
     $J_{i,Y}$ which can cover a given $j$.
     Furthermore, this argument shows that in any hierarchy with
     initial marking $\mu_{i}$, then either 
     $d_{\C(Y)}(\mu_{i},\mu_{i,j})<2M$ for all $j$ or  
     the projection to $\C(Z)$ of the terminal marking of 
     this hierarchy is within $M$ of the point $\phihat(\mu_{i})$; 
     when the later condition is satisfied we say that this hierarchy 
     has \emph{traversed  $Z$}.

     What we shall now show is that the hierarchy from $\mu_i$ to
     $\mu'_i$ can not possibly traverse all the large domains
     appearing in the hierarchy $G_i$.  This then implies that there
     exists a large domain of $G_i$, which we will call $Y'_i$, in
     which $d_{\C(Y'_i)}(\mu,\mu'_i)<2M$.  The strongly bounded
     hypothesis implies that $\diam_{\C(Y'_i)}([x_i,y_i])<K$, together
     with the fact that $d_{\C(Y'_i)}(\mu,\phihat(\mu_i))>4M+D$ this 
     allows us to conclude that the main geodesic of any hierarchy
     from $\mu'_i$ to $[x_i,y_i]$ must pass through the vertex
     $\partial Y'_i\in\C(S)$.  From this and the fact that
     $[x_i,y_i]\subset F_0$ we then give a quick deduction of the theorem.

     Recall that $\G_{G_{i},4M+D}$ denotes the set of domains in
     $G_{i}$ supporting geodesics larger than $4M+D$; unless noted
     otherwise all time orderings on domains will refer to the
     ordering in the hierarchy $G_{i}$.  Recall that
     the set of large initial domains for $G_{i}$ has cardinality at
     most $2\xi(S)$.  Thus, by \fullref{lemmalargesetsexists} there
     exists $Y'_{i}\subsetneq S$ with $Y'_{i}\in \G_{G_{i},4M+D}$, such that
     if we define $$\mathcal{I}_{i}=\{Y'_{i}\}\cup \{Z\in\G_{G_{i},4M+D}:
     Y'_{i}\prec_{t} Z\},$$
     then for $\omega$--almost every $i$ we have
$\displaystyle\sum_{Z\in\mathcal{I}_{i}}d_{\C(Z)}(\mu_{i}, 
\phihat(\mu_{i}))>\alpha
d_{\M(S)}(\mu_{i},\phihat(\mu_{i}))$, for a fixed $\alpha>0$
independent of $i$.
     
     Since the $J_{i,Z}$ cover the path from $\mu_{i}$ to $\mu_{i}'$
     with degree at most $2\xi(S)$, we have
     $2k\xi(S)\geq\sum_{Z\in\mathcal{I}_i}|J_{i,Z}|$.  Now we use the 
     result from \cite[page 962]{MasurMinsky:complex2} that there exists a
     constant $\beta>0$ (depending only on the topological type of
     $S$) for which $|J_{i,Z}|\geq \beta
     d_{\C(Z)}(\mu_{i},\phihat(\mu_{i}))$ if the path from $\mu_{i}$
     to $\mu_{i}'$ traverses the large domain $Z$.  Thus, if all the
     domains $Z\in \mathcal{I}_i$ are traversed between $\mu_i$ and
     $\mu'_i$, we have $$2k\xi(S)\geq \sum_{Z\in\mathcal{I}_i}|J_{i,Z}|\geq
     \sum_{Z\in\mathcal{I}_i}\beta d_{\C(Z)}(\mu_{i},\phihat(\mu_{i}))
     \geq \alpha\beta d_{\M(S)}(\mu_{i},\phihat(\mu_{i})).$$
     Thus, $k\geq \frac{\alpha\beta
     d_{\M(S)}(\mu_{i},\phihat(\mu_{i}))}{2\xi(S)}$; taking $c\leq
     \frac{\alpha\beta}{2\xi(S)}$
     we obtain a contradiction, and hence we may assume that not all
     the large domains in $\mathcal{I}_i$ are traversed between $\mu_i$
     and $\mu'_i$.
     
     In particular, we may assume $Y'_{i}$ can not be completely
     traversed between $\mu_i$ and $\mu'_i$ and is thus a large domain for
     any hierarchy from $\mu_{i}'$ to $\phihat(\mu_{i})$.  Therefore,
     the main geodesic of any hierarchy from $\mu_{i}'$ to
     $\phihat(\mu_{i})$ contains the vertex $\partial Y'_{i}$.  If
     $d_{\C(S)}(\partial Y'_{i}, \phihat(\mu_{i}))>\delta$, then any
     $\C(S)$ geodesic from $\pi_{\C(S)}(\mu_{i}')$ to the main
     geodesic of $H_{i}$ must pass within $\delta$ of $\partial Y'_{i}$.
     This yields a uniform bound on $\diam(C_{\mu_{i}}\cup
     C_{\mu_{i}'})$ and the same argument given before shows that in
     this case $\Phi(\mu)=\Phi(\mu')$.
     
     Thus, we may assume that the only choices for $Y'_{i}$ are those 
     which satisfy\break $d_{\C(S)}(\partial Y'_{i}, \phihat(\mu_{i}))\leq
     \delta$. It then follows from the Move distance and projections 
     Theorem that  
     for $\omega$--almost every $i$ there exists a domain $Q_{i}$ for which
     $d_{\C(S)}(\partial Q_{i}, \phihat(\mu_{i}))\leq \delta$ and with
     $\ulim d_{\C(Q_{i})}(\mu_{i},\phihat(\mu_{i}))\to \infty$. 
      On the other hand, the strong boundedness  of the pair  $(x,y)$ 
     implies that 
     $\diam_{\C(Q_{i})}([x_i,y_i])<D$ for $\omega$--almost every $i$, 
     hence  by \fullref{MMtwLLL}, 
     any hierarchy between $\mu_{i}'$ and any  point on $[x_{i},y_{i}]$ 
     must have a geodesic supported on $Q_{i}$.
     Thus the main geodesic of any hierarchy between $\mu_{i}'$ and
     any point on $[x_{i},y_{i}]$ must pass through the point
     $\partial Q_{i}$.  This provides a uniform bound on the distance
     between $\phihat(\mu_{i})$ and $\phihat(\mu_{i}')$ and therefore
     shows that $\Phi(\mu)=\Phi(\mu')$, completing the proof.

(3)\qua  Any bi-Lipschitz flat which intersects $F_{0}$ in more 
     than one point
     contains an arc of $F_{0}$. If this flat has rank greater
     than one, then this arc must be part of an embedded circle. But
     this would contradict the fact that there is a unique path in
$\coneM$
     between each pair of points of $F_{0}$.

     Thus a bi-Lipschitz flat of rank greater than one can
     intersect $F_{0}$ in at most one point.
\end{proof}

\section{$\cone(\MCG)$ is tree-graded}\label{sectiontreegraded}

We begin by stating 
a consequence of \fullref{Fzertree} and use this as a segue to 
mention some interesting questions about the mapping class group.
In \fullref{Fzertree} we obtain a unique path in $\coneM$ between
any pair of points in $F_{0}$. From this it follows that the removal 
of any point of $F_{0}$ disconnects $\coneM$, ie, 
if $S$ is any surface and $\coneM(S)$ is any asymptotic cone of 
$\MCG(S)$, then  $\coneM(S)$ has a global cut-point. 
Since asymptotic cones of groups are homogeneous, this can be
stated as:
\begin{thm}\label{cutpointstw}If $S$ is any surface and 
    $\coneM(S)$ is any asymptotic cone of 
    $\MCG(S)$, then every point of  $\coneM(S)$ is 
    a global cut-point.
\end{thm}

We now mention the following definition which was introduced in
\cite{DrutuSapir:TreeGraded}.

\begin{defn} Let ${\mathbb F}$ be a complete geodesic metric space and
    let $\mathcal{P}$ be a collection of closed geodesic subsets
    (called \emph{pieces}). 
    The space  ${\mathbb F}$ is said to be \emph{tree-graded 
    with respect to  $\mathcal{P}$} 
    when the following two properties are satisfied:
    \begin{itemize}
        \item The intersection of each pair of pieces has at most one 
	point.
        \item Every simple geodesic triangle in ${\mathbb F}$ is
	contained in one piece.
    \end{itemize}  
\end{defn}

In Remark~2.31 
of \cite{DrutuSapir:TreeGraded}, it is noted ``every 
asymptotic cone of a group which has a cut-point is tree graded with
respect to a uniquely determined collection of pieces each of which
is 
either a singleton or a closed geodesic subset 
without cut-points.'' In light of this remark,
\fullref{cutpointstw} implies: 
\begin{cor}\label{MCGtreegraded} For every surface $S$, 
       each asymptotic cone, $\coneM(S)$, is tree-graded.
\end{cor}
\noindent This result 
raises the question of describing the pieces. A useful object in the 
study of tree-graded structures is the \emph{transversal tree} 
containing a point $\mu$: this is defined to be the maximal tree 
(possibly consisting of only the point $\mu$), 
satisfying the property that it intersects each piece in at most one 
point. One consequence of 
\fullref{Fzertree} is that $F_{0}$ is a subset of the transversal 
tree containing the basepoint $0$. (Similarly by defining $F_{\mu}$ 
for any $\mu\in\coneM(S)$, we obtain a subset of the transversal tree 
containing the point $\mu$.) 

One of the questions solved by \cite{DrutuSapir:TreeGraded} 
is a proof of quasi-isometric rigidity for groups whose 
asymptotic cones admit tree-gradings of a certain form, namely the
asymptotic cones of their peripheral subgroups are cut-point free 
(see \cite{BehrstockDrutuMosher:thick} for a generalization). 
This approach might shed light on the
following, which remains a major outstanding conjecture in the
study of the mapping class group:
\begin{conj} \label{rigid}

     Let $S$ be a non-exceptional surface of finite type. For any 
    finitely generated group $G$ quasi-isometric to the mapping class 
    group $\MCG(S)$, there exists a homomorphism 
    $G\to\MCG(S)$ with finite kernel and finite index image.
\end{conj}

\noindent

Informally, this question asks if the mapping class group 
is uniquely determined among finitely generated groups by its large 
scale geometry. (See \cite{Mosher:Rigid1punct} for a further 
discussion and a proof of \fullref{rigid} for mapping class 
groups of once punctured surfaces.)   

Another interesting  question is to determine: to
what extent are mapping class groups relatively hyperbolic.
There are several competing definitions for  relative hyperbolicity,
one of which is that the Cayley graph ``electrified'' over the cosets 
of a finite collection of subgroups is
$\delta$--hyperbolic, this is now often called \emph{weak relative
hyperbolicity}. Another is weak relative hyperbolicity combined with a
hypothesis called the  Bounded Coset Penetration Property (BCP), which
places restrictions on how paths can travel between these cosets;
this 
condition was first formulated by Farb \cite{Farb:RelHyp} and is
often called \emph{(strong) relative hyperbolicity}. Many people have 
since studied another equivalent versions of relative hyperbolicity 
which was formulated by Gromov \cite{Gromov:hyperbolic}
and then elaborated on by Bowditch, Osin, Szczepa\'{n}ski, and others 
\cite{Bowditch:RelHyp,Osin:RelHyp,Szczepanski:relhyp}.

Another main result in \cite{DrutuSapir:TreeGraded} 
is that a group, $G$, is relatively hyperbolic with respect to a
collection of subgroups, $H_{1},\ldots,H_{n}$, which satisfy BCP 
if and only if 
$\cone(G)$  is tree-graded with respect to ultralimits of cosets of
the subgroups $H_{i}$.
They ask the open question of whether any group whose asymptotic cones
have
cut-points is relatively hyperbolic (See Problem~1.18 of 
\cite{DrutuSapir:TreeGraded}).
Masur and Minsky showed that 
$\MCG(S)$ is weakly hyperbolic relative to
stabilizers of multicurves \cite{MasurMinsky:complex1}, thus 
one knows that the maximal subsets without
cut-points in $\coneM(S)$ are not ultralimits of these subgroups. It 
is easy to verify that multicurve stabilizers do not satisfy BCP, 
nonetheless, there could be another collection of subgroups with 
respect to which the mapping class group is relatively hyperbolic. 
Accordingly, in an earlier version of this paper we asked the 
fundamental question: 
\begin{qn}\label{BCPquestion} Is there a collection of 
    subgroups of $\MCG(S)$ with
    respect to which this group is strongly relatively hyperbolic?
\end{qn}

An answer to this question can now be given. The first author in
joint work with Dru\c{t}u and Mosher
\cite{BehrstockDrutuMosher:thick} and independently work of Anderson,
Aramayona, and
Shackleton \cite{AASh:RelHypMCG} gives a complete
answer to this question by showing that for
no family of subgroups is the mapping class group strongly relatively 
hyperbolic.

\fullref{MCGtreegraded} contrasts with
the result of \cite{BehrstockDrutuMosher:thick} and
\cite{AASh:RelHypMCG} 
 to yield an example of a group whose 
asymptotic cone is tree-graded, but the group is not strongly
relatively hyperbolic
for any choice of subgroups; this answers Problem~1.18 of 
\cite{DrutuSapir:TreeGraded} and does so with a finitely presented
group. In fact, in \cite{BehrstockDrutuMosher:thick} a stronger result is 
shown, namely, there is no collection of subsets of the 
mapping class group for
which $\cone(\MCG)$ is tree-graded with  pieces coming from
ultralimits of these subsets; thus the result that  $\cone(\MCG)$ is, 
nonetheless, tree-graded answers another question of
\cite{DrutuSapir:TreeGraded} 
concerning whether pieces in the asymptotic cone always arise as
ultralimits of subsets.

We note the following question which this raises.
\begin{qn} What other groups besides $\MCG(S)$ have tree-graded
asymptotic cones, but are not strongly relatively hyperbolic with respect to 
   any subgroups? \ldots subsets?
\end{qn}

After learning of our methods, B.~Kleiner informed us that using
analogous
methods one can show that graph manifolds satisfy the property 
of having tree-graded asymptotic cones, but are not strongly
relatively 
hyperbolic with respect to any subgroups \cite{Kleiner:personal}.
(See \cite{KapovichKleinerLeeb:QIdeRham} for a related discussion. See also 
\cite{BehrstockDrutuMosher:thick} for a proof of the latter fact.)

The results about asymptotic cones throughout this work have all been 
independent of the choice of ultrafilter. In \cite{KSTT} it is shown
that if one assumes the negation of the 
Continuum Hypothesis, then there are uniform
lattices in certain semi-simple Lie groups with $2^{2^{\omega}}$
different asymptotic cones up to homeomorphism (here $\omega$ denotes 
the first uncountable ordinal). Even in the presence
of the Continuum Hypothesis, \cite{DrutuSapir:TreeGraded} produce
examples of finitely generated groups with $2^{\omega}$
non-homeomorphic asymptotic cones. These results make natural the
following two questions.
\begin{qn}\label{qnmanycones}
    \noindent
    \begin{enumerate}
    \item How many non-homeomorphic asymptotic cones does the mapping
class group 
    have?
    \item (Sapir)\qua How many non-bi-Lipschitz equivalent asymptotic
cones does
    the mapping class group have?
\end{enumerate}
\end{qn}

Note that for the examples given in \cite{DrutuSapir:TreeGraded} of
finitely generated groups with $2^{\omega}$ non-homeomorphic
asymptotic cones, the fundamental group can be used to distinguish the
different asymptotic cones.  Producing non-homeomorphic asymptotic
cones for the mapping class group would require more subtle techniques
since automaticity of the mapping class group implies that every
asymptotic cone for this group is simply connected---this follows from
results in \cite{Papasoglu:ConesQuadratic} and
\cite{Mosher:automatic}.

\bibliographystyle{gtart}
\bibliography{link}

\end{document}